\documentclass[aap,MSNbibl,seceqn,citesort,dvips]{arximspdf}
\usepackage{mathbh}
\usepackage{mathrsfs}

\doi{10.1214/12-AAP849} %kopijuoti is PTS
\volume{23}
\issue{2}
\pubyear{2013}
\firstpage{637}
\lastpage{664}

\makeatletter
\newcommand{\1}{\mathbh{1}}
\newcommand{\B}{\mathscr{B}}
\newcommand{\da}{{\downarrow}}
\newcommand{\lra}{\longrightarrow}
\newcommand{\lmt}{\longmapsto}
\newcommand{\bsb}{\bolds}
\newcommand{\bsm}{\mathbf}
\renewcommand{\P}{\mathbb{P}}
\newcommand{\E}{\mathbb{E}}
\newcommand{\mc}{\mathcal}
\newcommand{\dxi}{{\widehat{\xi}}}
\newtheorem{thmi}{Theorem}
\newtheorem{cori}{Corollary}
\newtheorem{thmm}{Theorem}[section]
\newtheorem{prop}[thmm]{Proposition}
\newtheorem{lem}[thmm]{Lemma}

\newproclaim{defi}[thmm]{Definition}
\newproclaim{prob}[thmm]{Problem}
\newproclaim{disc}[thmm]{Discussion}
\newproclaim{nota}{Notation}
\newproclaim{conj}{Conjecture}
\newproclaim{rmk}[thmm]{Remark}
\newproclaim{eg}[thmm]{Example}
\makeatother

\begin{document}
\begin{frontmatter}

\title{Sharp benefit-to-cost rules for the evolution of cooperation on regular graphs}
\runtitle{Sharpness of two simple rules}

\begin{aug}
\author{\fnms{Yu-Ting} \snm{Chen}\corref{}\ead[label=e1]{ytchen@math.ubc.ca}}
\runauthor{Y.-T. Chen}
\affiliation{University of British Columbia}
\address{Department of Mathematics\\
University of British Columbia\\
1984 Mathematics Rd.\\
Vancouver, British Columbia\\
Canada V6T 1Z2\\
\printead{e1}}
%adresu isvedimo komanda gale!
\end{aug}

% HISTORY:
\received{\smonth{7} \syear{2011}}
\revised{\smonth{1} \syear{2012}}

% ABSTRACT
%
\begin{abstract}
We study two of the simple rules on finite graphs under the death--birth
updating and the imitation updating discovered by Ohtsuki, Hauert,
Lieberman and Nowak [\textit{Nature} \textbf{441} (2006) 502--505]. Each rule
specifies a payoff-ratio cutoff point for the magnitude of fixation
probabilities of the underlying
evolutionary game between cooperators and defectors.
We view the Markov chains associated with the two updating mechanisms
as voter model perturbations. Then we present a first-order
approximation for fixation probabilities of general voter model
perturbations on finite graphs subject to small perturbation %weak
%selection
in terms of the voter model fixation probabilities.
In the context of regular graphs, we obtain algebraically explicit
first-order approximations for the fixation probabilities of
cooperators distributed as certain uniform distributions. These
approximations lead
to a rigorous proof that both of the rules of Ohtsuki et al. are valid
and are sharp.
\end{abstract}

% KEYWORDS
%


\end{frontmatter}

%s1 #&#
\section{Introduction}\label{secintro}
The main objective of this paper is to investigate the simple rules for
the evolution of cooperation by clever, but nonrigorous, arguments of
pair-approximation on certain large graphs in Ohtsuki, Hauert,
Lieberman and Nowak~\cite{OHLN}.
For convenience, we name these rules and their relatives
comprehensively as \textit{benefit-to-cost rules} (abbreviated as \textit{$ b / c$ rules})
for reasons which will become clear later on.
The work~\cite{OHLN} takes spatial structure into consideration and
gives an explanation with analytical criteria for the ubiquity of
cooperative entities observed in biological systems and human
societies. (See also the references in~\cite{OHLN} for other models on
structured populations.) In particular, this provides a way to overcome
one of the major difficulties in theoretical biology since Darwin. (See
Hamilton~\cite{HGESB}, Axelrod and Hamilton~\cite{AHEC}, Chapter 13
in Maynard Smith~\cite{MSETG} and many others.)

We start by describing the evolutionary games defined in~\cite{OHLN}
and set some definitions.
Consider a finite, connected, and simple (i.e., undirected and without
loops or parallel edges) graph $G=(\mathsf V,\mathsf E)$ on $N$
vertices. (See, e.g.,~\cite{BMGT} for the standard terminology of
graph theory.) Imagine the graph as a social network where a population
of $N$ individuals occupy the vertices of $G$ and the edges denote the
links between the individuals. The population consists of cooperators
and defectors labeled by $1$'s and $0$'s, respectively.
Their fitness is described through payoffs from encounters as follows.
Consider a $2\times2$ payoff matrix
%
%e1 #&#
\begin{equation}\label{eqPi}
\bsb\Pi=\pmatrix{
\Pi_{11} & \Pi_{10} \vspace*{2pt}\cr
\Pi_{01} & \Pi_{00} }
=
\pmatrix{
b-c & -c \vspace*{2pt}\cr
b & 0 }.
\end{equation}
Here, while positive constants are natural candidates for both
benefit $b$ and cost $c$,
we allow arbitrary reals for their possible values unless otherwise mentioned.
Each entry $\Pi_{ij}$ of $\bsb\Pi$ denotes the payoff that an
$i$-player receives from a $j$-player.
Hence, the payoff of a cooperator is $bn-ck$ if $n$ of its $k$
neighbors are cooperators, and the payoff of a defector is $bm$ if $m$
of its neighbors are cooperators.
Then the fitness of an individual located at $x$ is given by a convex
combination of the \textit{baseline fitness} with weight $1-w$ and its
payoff with weight $w$. The baseline fitness is normalized to $1$ for
convenience. The parameter $w$ is interpreted as the \textit{intensity of
selection}.
Therefore, \textit{weak selection} means that payoff is a small proportion
of fitness compared to the baseline fitness.

In contrast to game theory where strategies are decided by rational
players, evolutionary game theory considers the random evolution of
interacting players in which the ``fitter'' strategies have better
chances to replicate.
We study two of the updating mechanisms under weak selection in \cite
{OHLN} for the evolution of cooperation throughout our work.
Under the death--birth updating, we kill a random individual and then
let its neighbors compete for the vacant vertex
with success probability proportional to the fitness of its neighbors.
Under the imitation updating, a random individual updates its strategy,
but now it will either adhere to its original strategy or imitate one
of the neighbors' strategies with success probability proportional to fitness.
In this way, each updating mechanism defines a Markov chain on
configurations consisting of $1$'s and $0$'s, or more specifically, a
spin system in the sense of Liggett~\cite{LIPS} where each vertex can
adopt only two possible opinions, $1$ and $0$. Despite the simplicity
of the transition rates, the readers may observe that the spin systems
pose certain difficulty in terms of the classical approaches in
interacting particle systems. For example, as a result of the asymmetry
of payoffs, there is no symmetry between $1$'s and $0$'s in the two
spin systems. In addition,
it is not hard to see that in general the two spin systems are not
attractive.

We are now ready to describe the $b/c$ rules for the two evolutionary
games which are surprisingly simple criterions for criticality in the
asymptotic.
The degree of a graph is defined in~\cite{OHLN} to be the average
number of neighbors per vertex. Put a single cooperative mutant on the
vertices with a random location.
Then the clever, but nonrigorous, calculations in the supplementary
notes of~\cite{OHLN}, supported by several numerical
simulations, lead to the following \textit{b/c rule} for
the death--birth updating under weak selection on certain large graphs
of degree~$k$;
selection favors cooperation whenever $b/c>k$ and selection opposes
cooperation whenever $b/c<k$. Here, selection \textit{favors} (resp.,
\textit{opposes}) cooperation if the probability that a single cooperative
mutant converts the defecting population completely into a cooperative
population is strictly higher (resp., lower) than
the fixation probability $1/N$ of a neutral mutant. [See (\ref
{eqpiREG}) for the latter probability. Equation (\ref{eqpiREG}) also
shows, in particular, that if the graph is regular, then the fixation
probability of a neutral mutant at an arbitrary location without
further randomization is precisely $1/N$.] A similar $b/c$ rule under
the imitation mechanism is discussed in the supplementary notes
of~\cite{OHLN}, with the modification that the cutoff point $k$ should
be replaced by $k+2$.
We remark that the work~\cite{OHLN} also considers the birth--death
updating (in contrast to the death--birth updating) and its associated
$b/c$ rule. See~\cite{TDWEC} for a further study of these $b/c$ rules.
For more $b/c$ rules see~\cite{HGESB,NR} and~\cite{S}, to
name but a few. The monograph~\cite{NED} gives an authoritative and
excellent introduction to evolutionary dynamics.

Lying at the heart of the work~\cite{OHLN} to obtain selective
advantage of cooperators is the introduction of structured populations.
This is manifested by the role of a fixed degree as population size
becomes large. Consider instead a naive model where only fractions of
players in a large population are concerned and the same payoff matrix
(\ref{eqPi}) is in effect for evolutionary fitness. The fractions
$z_C$ and $z_D$ of cooperators and defectors are modeled through
replicator equations as
%
%e2 #&#
%e3 #&#
\begin{eqnarray}\label{eqnaive}
\dot{z}_C&=&z_C(\rho_C -\bar{\rho}),
\nonumber
\\[-8pt]
\\[-8pt]
\nonumber
\dot{z}_D&=&z_D(\rho_D-\bar{\rho}).
\end{eqnarray}
Here, by the equality $z_C+z_D=1$, the payoffs for cooperators and
defectors are $\rho_C=bz_C-c$ and $\rho_D=bz_C$, and $\bar{\rho}$
is the average payoff given by $z_C\rho_C+z_D\rho_D$. By (\ref
{eqnaive}), the fraction of cooperators satisfies the following
logistic differential equation:
\[
\dot{z}_C=-cz_C(1-z_C).
\]
Hence, any proper fraction of cooperators must vanish eventually
whenever cost~$c$ is positive. See, for example, Chapter 7 in \cite
{HSEGPD}, Chapter 4 in~\cite{NED} or Section 3 in~\cite{HED} for
this model and more details. As discussed in more detail later on, a~similar result
holds in any unstructured population of finite size under the
death--birth updating.
Informally, a spatial structure, on one hand, promotes the formation of
cliques of cooperators which collectively have a selective advantage
and, on the other hand, reduces the exploitation of cooperators by defectors.

In~\cite{ON}, Ohtsuki and Nowak gave a rigorous proof of the $b/c$
rules on large cycles under weak selection, in particular for the two
updating mechanisms. The results in~\cite{ON} exploit the fact that on
cycles the fixation probabilities under each updating mechanism satisfy
a system of birth--death-process type difference equations, and exact
fixation probabilities can be derived accordingly.
It is easy to get the exact solvability of fixation probabilities by
the same approach on complete graphs, although on each \textit{fixed}
one of these graphs cooperators are always opposed under weak selection
for the death--birth updating. (See~\cite{OHLN} and Remark~\ref{rmkegfs1}(3). Note that the degree of a complete graph has the same
order as the number of vertices.)
It seems, however, harder to obtain fixation probabilities by extending
this approach beyond cycles and complete graphs.

In this work, we will view each of the two spin systems as a \textit{voter model
perturbation} on a (finite, connected and simple) graph \textit{of
arbitrary size}. Voter model perturbations are studied in Cox and
Perkins~\cite{CPLV} and further developed in generality in Cox,
Durrett and Perkins~\cite{CDPVMP} on
transient integer lattices $\Bbb Z^d$ for $d\geq3$.
On the infinite
lattices considered in~\cite{CDPVMP}
(often sharp) conditions, based on a related reaction diffusion equation,
were found to ensure the
coexistence of $1$'s and $0$'s, or to ensure that one
type drives the other out. In particular, a rigorous proof of the $b/c$
rule under the death--birth updating on these infinite graphs is
obtained in~\cite{CDPVMP}. Here, in the context of finite graphs,
the voter model perturbation associated with each of the spin systems
fixates in one of the the two absorbing states of
all $1$'s or all $0$'s, and we give a first-order approximation to the
fixation probabilities by expansion. In spite of the apparent
differences in the
settings, there are interesting links between our fixation probability
expansions and the reaction diffusion equation criteria in~\cite{CDPVMP}.

Let us now introduce voter model perturbations as a family of spin
systems. Denote by $x$ a vertex and by $\eta$ a configuration.
Let $c(x,\eta)$ be the flipping rate of the (nearest-neighbor) voter
model at $x$ given $\eta$. Hence, $c(x,\eta)$ is equal to the
probability of drawing a neighbor of $x$ with the opposite opinion to
that of $x$.
Then interpreted narrowly in the context considered in this paper, the
rates of a voter model perturbation are given by
%
%e4 #&#
\begin{equation}
c^w(x,\eta)=c(x,\eta)+wh_{1-\eta(x)}(x,\eta)+w^2g_w(x,\eta)\geq
0\label{eqpf}
\end{equation}
for a small perturbation rate $w> 0$. Here, $h_1$, $h_0$ and $g_w$ for
all $w$ small are uniformly bounded. We refer the readers to Chapter V
in~\cite{LIPS} here and in the following for the classical results of
voter models and to Section 1 in~\cite{CDPVMP} for a general
definition of voter model perturbations on transient integer lattices.

We discuss in more detail the aforementioned result in~\cite{CDPVMP}
which is closely related to the present work.
A key result in~\cite{CDPVMP} states
that on the integer lattices $\Bbb Z^d$ for $d\geq3$, the invariant
distributions of a voter model perturbation subject to small
perturbation can be determined by the reaction function through a
reaction-diffusion PDE. (See Section 1.2 in~\cite{CDPVMP}.)
Here, the reaction function takes the form
%
%e5 #&#
\begin{equation}
u\lmt\lim_{s\to\infty}\int D(0,\eta)\P_{\mu_u}(\xi_s\in d\eta
), \qquad u\in[0,1],\label{eqrp}
\end{equation}
where $((\xi_s),\P_{\mu_u})$ denotes the voter model starting at the
Bernoulli product measure $\mu_u$ with density $\mu_u(\eta(x)=1)=u$,
and the difference kernel $D$ is defined by
%
%e6 #&#
\begin{equation}
D(x,\eta)=\widehat{\eta}(x)h_1(x,\eta)-\eta(x)h_0(x,\eta)\label{eqD}
\end{equation}
with $\widehat{\eta}(x)\equiv1-\eta(x)$.
By the duality between voter models and coalescing random walks, the
reaction function defined by (\ref{eqrp}) can be expressed explicitly
as a polynomial with coefficients consisting of coalescing
probabilities of random walks.

The justification in~\cite{CDPVMP} of the $b/c$ rule for the
death--birth updating on transient integer lattices is under a slightly
different definition, for the sake of adaptation to the context of
infinite graphs. Precisely, the result in~\cite{CDPVMP} states that
whenever $b/c>k$ (resp., $b/c<k$) and there is weak selection, given
infinitely many cooperators (resp., defectors) at the beginning, any
given finite set of vertices will become occupied by cooperators
(resp., defectors) from some time onward almost surely. Here, $k$
refers to the degree of each vertex in the underlying lattice and is
equal to~$2d$ on $\Bbb Z^d$.
The $b/c$ rule under the imitation updating on the same integer
lattices is verified in~\cite{CID}, under the same definition in
\cite{CDPVMP} except that, as pointed out in~\cite{OHLN}, the cutoff
point $k$ needs to be replaced by $k+2$.

We now discuss our result for voter model perturbations on (finite,
connected and simple) graphs of arbitrary size. We first work with
discrete-time Markov chains of the voter model perturbations. We
assume, in addition, the chain starting at any arbitrary state is
eventually trapped at either of the two absorbing states, the all-$1$
configuration $\mathbf1$ and the all-$0$ configuration $\mathbf0$.
This is a property enjoyed by both updating mechanisms under weak selection.
We proceed analytically and decompose the transition kernel $P^w$ of a
voter model perturbation with perturbation rate $w$ as the sum of the
transition kernel $P$ of the voter model and a signed kernel~$K^w$. We
apply an elementary expansion to finite-step transitions of the
$P^w$-chain in the spirit of Mayer's cluster expansion~\cite{MMMD} in
statistical mechanics. Then we show every linear combination of the
fixation probabilities of $\mathbf1$ and $\mathbf0$ subject to small
perturbation admits an infinite series expansion closely related to the
voter model.
A slight refinement of this expansion leads to our main result for the
general voter model perturbations on which our study of the $b/c$ rules relies.

Precisely, our main result for the voter model perturbations on finite
graphs (Theorem~\ref{thmTaylor}) can be stated as follows. Regard the
voter model perturbation with perturbation rate $w$ as a
continuous-time chain $(\xi_s)$ with rates given by (\ref{eqpf}).
Recall by our assumption that the chain starting at any state is
eventually trapped at either of the absorbing states $\mathbf1$ and
$\mathbf0$.
Define $\tau_{\mathbf1}$ for the time to the absorbing state~$\mathbf
1$ and
%
%e7 #&#
\begin{equation}
\overline{H}(\xi)=\sum_{x\in\mathsf V}H(x,\xi)\pi(x)\label{eqHbar}
\end{equation}
for any $H(x,\xi)$. Here, $\pi(x)$ is the invariant distribution of
the (nearest-neighbor) random walk on $G$ given by
%
%e8 #&#
\begin{equation}
\pi(x)=\frac{d (x)}{2\cdot\#\mathsf E}\label{eqinv}
\end{equation}
with $d (x)$ being the degree of $x$, that is, the number of neighbors
of $x$. (See, e.g.,~\cite{AFIPS} and~\cite{LRWG} for random walks
on graphs.)
Then as $w\lra0+$, we have the following approximation:
%
%e9 #&#
\begin{equation}
\P^w(\tau_{\mathbf1}<\infty)=\P(\tau_{\mathbf
1}<\infty)+w\int_0^\infty\Bbb E[\overline{D}(
\xi_s)]\,ds+O(w^2).\label{eqapprox}
\end{equation}
Here, $\P^w$ and $\P$ (with expectation $\E$) denote the laws of the
voter model perturbation with perturbation rate $w$ and the voter
model, respectively, both subject to the same, but arbitrary, initial
distribution, and $D$ is the difference kernel defined by (\ref
{eqD}). Moreover, the integral term on the right-hand side of (\ref
{eqapprox}) makes sense because
$\overline{D}(\xi_s)\in L_1(d\Bbb P\otimes ds)$.

We apply the first-order approximation (\ref{eqapprox}) to the two
evolutionary games only on regular graphs. (A graph is a $k$-regular if
all vertices have the same degree $k$ and a graph is regular if it is
$k$-regular graph for some $k$.)
Under weak selection, the approximation (\ref{eqapprox}) implies that
we can approximate $\P^w(\tau_{\mathbf1}<\infty)$ by $\P(\tau
_{\mathbf1}<\infty)$ and the $0$-potential of $\overline{D}$
%
%e10 #&#
\begin{equation}
\int_0^\infty\Bbb E[\overline{D}(\xi_s)]\,ds,%=\int_0^
\label{eqint}
\end{equation}
all subject to the same initial distribution. Moreover, the comparison
of $\P^w(\tau_{\mathbf1}<\infty)$ for small $w$ with $\P(\tau
_{\mathbf1}<\infty)$ is possible whenever the $0$-potential
is nonzero, with the order determined in the obvious way.
For notions to be introduced later on, we take as initial distribution
the uniform distribution $\mathbf u_n$ on the set of configurations
with exactly $n$ many $1$'s, where $1\leq n\leq N-1$.
Each $0$-potential in~(\ref{eqint}) starting at $\mathbf u_n$
can be derived from the same $0$-potentials starting at Bernoulli
product measures $\mu_u$ with density $u\in[0,1]$. Furthermore, each
$0$-potential with starting measure $\mu_u$ can be expressed in terms
of some (expected) coalescing times of coalescing random walks. This is
in contrast to the involvement of coalescing probabilities for the
reaction functions in the context in~\cite{CDPVMP}. By resorting to a
simple identity in~\cite{AFIPS} between meeting times and hitting
times of random walks, we obtain explicit forms of the coalescing times
involved. Hence, by (\ref{eqapprox}), we obtain the fixation
probabilities explicit up to the first-order term, and the precise
result is stated as follows.

%th1 #&#
\begin{thmi}\label{thmi1}
Let $G$ be any (finite, connected and simple) graph on $N$ vertices.
Suppose in addition that $G$ is $k$-regular, that is, every vertex of
$G$ has precisely degree $k$. Fix $1\leq n\leq N-1$ and let $\P
^w_{\mathbf u_n}$ denote the law of the particular evolutionary game
with small intensity of selection $w$ and initial distribution $\mathbf u_n$.
\begin{longlist}[(1)]
\item[(1)] Under the death--birth updating, we have
\[
\P^w_{\mathbf u_n}(\tau_{\mathbf1}<\infty)=\frac{n}{N}+ w
\biggl[\frac{kn(N-n)}{2N(N-1)}\biggr]
\biggl[\biggl(\frac{b}{k}-c\biggr)(N-2)+b\biggl(\frac
{2}{k}-2\biggr)\biggr]+O(w^2)
\]
as $w\lra0+$.
\item[(2)] Under the imitation updating, we have
\begin{eqnarray*}
\P^w_{\mathbf u_n}(\tau_{\mathbf1}<\infty)&=&\frac{n}{N}+w
\biggl[\frac{k(k+2)n(N-n)}{2(k+1)N(N-1)}\biggr]\\
&&\hspace*{22pt}{}\times
\biggl[\biggl(\frac{b}{(k+2)}-c\biggr)(N-1)-\frac{
(2k+1)b-ck}{k+2}\biggr]+O(w^2)
\end{eqnarray*}
as $w\lra0+$.
\end{longlist}
Here, in either $(1)$ or $(2)$, we use Landau's notation $O(w^2)$ for a
function $\theta(w)$ such that $|\theta(w)|\leq Cw^2$ for all small
$w$, for some constant $C$ depending only on the graph $G$ and the
particular updating mechanism.
\end{thmi}

Before interpreting the result of Theorem~\ref{thmi1} in terms of
evolutionary games, we first introduce the following definition which
is stronger than that in~\cite{OHLN}. We say selection \textit{strongly}
favors (resp., opposes) cooperation if for \textit{every} nontrivial
$n$, that is, $1\leq n\leq N-1$, the following holds: The probability
that $n$ cooperative mutants with a joint location distributed as
$\mathbf u_n$ converts the defecting population completely into a
cooperative population is strictly higher (resp., lower) than $n/N$.
[Here, $n/N$ is the fixation of probability of $n$ neutral mutants
again by (\ref{eqpiREG}).] Under this definition, Theorem~\ref{thmi1} yields
simple algebraic criteria for both evolutionary games stated as follows.

%co1 #&#
\begin{cori}\label{cori1}
Suppose again that the underlying social network is a $k$-regular graph
on $N$ vertices.
\begin{longlist}[(1)]
\item[(1)] For the death--birth updating, if
\[
\biggl(\frac{b}{k}-c\biggr)(N-2)+b\biggl(\frac
{2}{k}-2\biggr)>0\qquad  (\mbox{resp., }<0),
\]
then selection \textit{strongly} favors (resp., opposes) cooperation
under weak selection.

\item[(2)] For the imitation updating, if
\[
\biggl(\frac{b}{(k+2)}-c\biggr)(N-1)-\frac{
(2k+1)b-ck}{k+2}>0 \qquad(\mbox{resp., }<0),
\]
then selection \textit{strongly} favors (resp., opposes) cooperation
under weak selection.
\end{longlist}
\end{cori}

Applied to cycles, the algebraic criteria in Corollary~\ref{cori1}
under the aforementioned stronger definition coincide with the
algebraic criteria in~\cite{ON} for the respective updating mechanism.
See also equation (3) in~\cite{TDWEC} for the death--birth updating.

As an immediate consequence of Corollary~\ref{cori1}, we have the
following result.

%co2 #&#
\begin{cori}\label{cori2}
Fix a degree $k$.
\begin{longlist}[(1)]
\item[(1)] Consider the death--birth updating. For every fixed pair $(b,c)$
satisfying $b/k>c$ (resp., $b/k<c$), there exists a positive integer
$N_0$ such that on any $k$-regular graph $G=(\mathsf V,\mathsf E)$ with
$\#\mathsf V\geq N_0$, selection \textit{strongly} favors (resp.,
opposes) cooperation under weak selection.

\item[(2)] Consider the imitation updating. For every fixed pair $(b,c)$
satisfying $b/(k+2)>c$ [resp., $b/(k+2)<c$], there exists a positive
integer $N_0$ such that on any $k$-regular graph $G=(\mathsf V,\mathsf
E)$ with $\#\mathsf V\geq N_0$, selection \textit{strongly} favors
(resp., opposes) cooperation under weak selection.
\end{longlist}
\end{cori}

In this way, we rigorously prove
the validity of the $b/c$ rule in~\cite{OHLN} under each updating
mechanism. It is in fact a \textit{universal} rule valid for any
nontrivial number of cooperative mutants, and holds \textit{uniformly}
in the number of vertices, for large regular graphs with a fixed degree
under weak selection.

%re1 #&#
\begin{rmk}\label{rmkegfs1}
(1) Although we only consider payoff matrices of the special form (\ref
{eqPi}) in our work, interests in evolutionary game theory do cover
general $2\times2$ payoff matrices with arbitrary entries. (See, e.g.,
\cite{OHLN} and~\cite{ON}.) In this case, a general $2\times2$
matrix $\bsb\Pi^*=(\Pi^*_{ij})_{i,j=1,0}$ is taken to define payoffs
of players with an obvious adaptation of payoffs under $\bsb\Pi$. For
example, the payoff of a cooperator is $(\Pi^*_{11}-\Pi^*_{10})n+k\Pi
^*_{10}$ if $n$ of its $k$ neighbors are cooperators. In particular, if
$\bsb\Pi^*$ satisfies\vadjust{\goodbreak}
the \textit{equal-gains-from-switching} condition (Nowak and Sigmund
\cite{NSESSPD})
%
%e11 #&#
\begin{equation}
\Pi_{11}^*-\Pi^*_{10}=\Pi^*_{01}-\Pi^*_{00},\label{eqegfs}
\end{equation}
then the results in Theorem~\ref{thmi1}, Corollaries~\ref{cori1} and~\ref{cori2} still hold
for $\bsb\Pi^*$ by taking $\bsb\Pi$ in their statements to be the
``adjusted'' payoff matrix
%
%e12 #&#
\begin{equation}\label{eqPiadj}
\bsb\Pi^a:=\pmatrix{
\Pi^*_{11}-\Pi^*_{00} & \Pi^*_{10}-\Pi^*_{00} \vspace*{2pt}\cr
\Pi^*_{01}-\Pi^*_{00} & 0
},
\end{equation}
which is of the form in (\ref{eqPi}).
See Remark~\ref{rmkegfs2} for this reduction.\vspace*{-6pt}

\begin{longlist}
\item[(2)] We stress that when $n=1$ or $N-1$ and the graphs are
vertex-transitive~\cite{BMGT} (and hence, regular) such as tori, the
exact locations of mutants become irrelevant. It follows that the
randomization by $\mathbf u_n$ is redundant in these cases.

\item[(3)]Let $G$ be the complete graph on $N$ vertices so that the
spatial structure is irrelevant.
Consider the death--birth updating and the ``natural case'' where benefit
$b$ and cost $c$ are both positive.
With the degree $k$ set equal to $N-1$, Theorem~\ref{thmi1}(1) gives
for any $1\leq n\leq N-1$ the approximation
\[
\P^w_{\mathbf u_n}(\tau_{\mathbf1}<\infty)=\frac
{n}{N}+w\frac{n(N-n)}{2N}\biggl[-c(N-2)-\biggl(2-\frac{N}{N-1}
\biggr)b\biggr]+O(w^2)
\]
as $w\lra0+$. Hence, cooperators are always opposed under weak
selection when $N\geq3$.
\end{longlist}
\end{rmk}

The paper is organized as follows.
In Section~\ref{secbs}, we set up the standing assumptions of voter
model perturbations considered throughout this paper and discuss their
basic properties. The Markov chains associated with the two updating
mechanisms in particular satisfy these standing assumptions, as stated
in Propositions~\ref{propdbu} and~\ref{propiu}.
In Section~\ref{secexp}, we continue to work on the general voter
model perturbations. We develop an expansion argument to obtain an
infinite series expansion of fixation probabilities under small
perturbation rates (Proposition~\ref{propexp}) and then refine its
argument to get the first-order approximation (\ref{eqapprox})
(Theorem~\ref{thmTaylor}). In Section~\ref{secbc}, we return to our
study of the two evolutionary games and give the proof of Theorem~\ref{thmi1}.
The vehicle for each explicit result is a simple identity between
meeting times and hitting times of random walks.
Finally, the proofs of Propositions~\ref{propdbu} and~\ref
{propiu} (that both updating mechanisms define voter model
perturbations satisfying our standing assumptions) are deferred to
Section~\ref{secuvmp}.

%s2 #&#
\section{Voter model perturbations}\label{secbs}
Recall that we consider only finite, connected and simple graphs in
this paper. Fix such a graph $G=(\mathsf V,\mathsf E)$ on $N=\# \mathsf
V$ vertices. Write $x\sim y$ if $x$ and $y$ are neighbors to each
other, that is, if there is an edge of $G$ between $x$ and $y$. We put
$d (x)$ for the number of neighbors of $x$.

Introduce an auxiliary number $\lambda\in(0,1]$.
Take a nearest-neighbor discrete-time voter model with transition probabilities
%
%e13 #&#
%e14 #&#
\begin{eqnarray}
\label{eqvoter}
P(\eta,\eta^x)&=&\frac{\lambda}{N}c(x,\eta), \qquad x\in\mathsf V,
\nonumber
\\[-8pt]
\\[-8pt]
\nonumber
P(\eta,\eta)&=&1-\frac{\lambda}{N}\sum_{x}c(x,\eta).
\end{eqnarray}
Here, $\eta^x$ is the configuration obtained from $\eta$ by giving up
the opinion of $\eta$ at $x$ for
\[
\widehat{\eta}(x):=1-\eta(x)
\]
and holding the opinions at other vertices fixed and we set
\[
c(x,\eta)=\frac{\#\{y\sim x;\eta(y)= \widehat{\eta}(x)\}}{d (x)}.
\]

We now define the discrete-time voter model perturbations considered
throughout this paper as follows. Suppose that we are given functions
$h_{i}$ and $g_w$ and a constant $w_0\in(0,1)$ satisfying
%
%e15 #&#
%e16 #&#
%e17 #&#
{\renewcommand{\theequation}{A.\arabic{equation}}
\setcounter{equation}{0}
\begin{eqnarray}
&\displaystyle\sup_{w\in[0,w_0],x,\eta}\bigl(|h_1(x,\eta)|+|h_0(x,\eta
)|+|g_w(x,\eta)|\bigr)\leq C_0<\infty.&
\label{A1}\\
&c^w(x,\eta):=c(x,\eta)+wh_{1-\eta(x)}(x,\eta)+w^2g_w(x,\eta)\geq0,&
\label{A2}\\
%&c^w(x,\eta)\leq1 \mbox{ for all $\eta$},\tag{A3}\label{A3}\\
&\hspace*{-53pt}c^w(x,\mathbf1)=c^w(x,\mathbf0)\equiv0 \qquad\mbox{for each $x\in
\mathsf V$}& %\tag{A3}
\label{A3}
\end{eqnarray}}
\hspace*{-2pt}for each $w\in[0,w_0]$.
Here, $\mathbf1$ and $\mathbf0$ denote the all-$1$ configuration and
the all-$0$ configuration, respectively.
In (\ref{A2}), we set up a basic perturbation of voter model rates up
to the second order. In terms of the voter model perturbations defined
below by $c^w(x,\eta)$,
we will be able to control the higher order terms in an expansion of
fixation probabilities with the uniform bound imposed in (\ref{A1}).
The assumption~(\ref{A3}) ensures that the voter model perturbations
have the same absorbing states $\mathbf1$ and $\mathbf0$ as the
previously defined voter model.

Under the assumptions (\ref{A1})--(\ref{A3}), we define for each
perturbation rate $w\in[0,w_0]$ a \textit{voter model perturbation} with
transition probabilities
%
%e18 #&#
%e19 #&#
\setcounter{equation}{1}
\begin{eqnarray}\label{eqcw}
P^w(\eta,\eta^x)&=&\frac{\lambda}{N}c^w(x,\eta),\qquad  x\in\mathsf
V,
\nonumber
\\[-8pt]
\\[-8pt]
\nonumber
P^w(\eta,\eta)&=&1-\sum_x\frac{\lambda}{N}c^w(x,\eta).
\end{eqnarray}
[Here we assume without loss of generality by (\ref{A1}) that each
$P^w(\eta,\cdot)$ is truly a probability measure, in part explaining
the need of the auxiliary number $\lambda$.]
In particular $P^0\equiv P$.

%pa2.0.0.1 #&#
\begin{nota*}
We shall write $\Bbb P^w_\nu$ for the law of the
voter model perturbation with perturbation rate $w$ and initial
distribution $\nu$ and set $\Bbb P_\nu:=\P_\nu^0$. In particular we
put $\Bbb P^w_{\eta}:=\P^w_{\delta_\eta}$ and $\Bbb P_\eta:=\P
_{\delta_\eta}$, where $\delta_\eta$ is the Dirac measure at $\eta$.
The discrete-time and continuous-time coordinate processes on $\{1,0\}
^{\mathsf V}$ are denoted by $(\xi_n;n\geq0)$ and $(\xi_s;s\geq0)$,
respectively. Here, and in what follows, we abuse notation to read
``$n$'' and other indices for the discrete time scale and ``$s$'' for
the continuous time scale whenever there is no risk of confusion.
\end{nota*}

Our last assumption, which is obviously satisfied by the $P$-voter
model thanks to the connectivity of $G$, is
%
%e20 #&#
{\renewcommand{\theequation}{A.\arabic{equation}}
\setcounter{equation}{3}
\begin{equation}\label{A4}
\P^w_\eta(\xi_n\in\{\mathbf1,\mathbf0\}\mbox{ for some
}n)>0\qquad \mbox{for every $\eta\in\{1,0\}^{\mathsf V}$}
\end{equation}}
\hspace*{-2pt}for each $w\in(0,w_0]$. Since $\mathbf1$ and $\mathbf0$ are
absorbing by the condition (\ref{A3}), it follows from the Markov
property that the condition (\ref{A4}) is equivalent to the condition that
the limiting state exists and can only be either of the absorbing
states $\mathbf1$ and $\mathbf0$
under $\Bbb P^w$ for any $w\in(0,w_0]$.

%pr2 #&#
\begin{prop}[(\cite{CDPVMP})]\label{propdbu}
Suppose that the graph is $k$-regular. Then
the Markov chain associated with the death--birth updating with small
intensity of selection $w$
is a voter model perturbation with perturbation rate $w$ satisfying
\textup{(\ref{A1})--(\ref{A4})} with $\lambda=1$ and
%
%e21 #&#
%e22 #&#
\setcounter{equation}{2}
\begin{eqnarray}
\label{eqdbh}
h_1&=&-(b+c)kf_0f_1+kbf_{00}+kf_0(bf_{11}-bf_{00}),
\nonumber
\\[-8pt]
\\[-8pt]
\nonumber
h_0&=&-h_1.
\end{eqnarray}
Here,
%
%e23 #&#
%e24 #&#
\begin{eqnarray}\label{eqf}
f_i(x,\eta)&=&\frac{1}{k}\#\{y;y\sim x, \eta(y)=i\},
\nonumber
\\[-8pt]
\\[-8pt]
\nonumber
f_{ij}(x,\eta)&=&\frac{1}{k^2}\#\{(y,z);x\sim y\sim z, \eta
(y)=i,\eta(z)=j\}.
\end{eqnarray}
\end{prop}
%
%
%pr3 #&#
\begin{prop}\label{propiu}
Suppose that the graph is $k$-regular. Then
the Markov chain associated with the imitation updating with small
intensity of selection $w$ is a voter model perturbation with
perturbation rate $w$ satisfying \textup{(\ref{A1})--(\ref{A4})} with
$\lambda=\frac{k}{k+1}$ and
%
%e25 #&#
%e26 #&#
\begin{eqnarray}
\label{eqiuh}
h_1&=&k[(b-c)f_{11}-cf_{10}]-\frac{k^2}{k+1}f_1[(b-c)f_{11}-cf_{10}+bf_{01}]\nonumber\\
&&{}-
\frac{k}{k+1}bf_1^2,
\nonumber
\\[-8pt]
\\[-8pt]
\nonumber
h_0&=&kbf_{01}-\frac{k^2}{k+1}f_0[(b-c)f_{11}-cf_{10}+bf_{01}]\\
&&{}-\frac{k}{k+1}
f_0[(b-c)f_1-cf_0],\nonumber
\end{eqnarray}
where $f_i$ and $f_{ij}$ are as in $(\ref{eqf})$.
\end{prop}

The proofs of Propositions~\ref{propdbu} and~\ref{propiu} are deferred to Section~\ref{secuvmp}.

\textit{The assumptions \textup{(\ref{A1})--(\ref{A4})} are in force from
now on.}

Let us consider some basic properties of the previously defined
discrete-time chains. First, as has been observed, we know that
\[
1=\P^w(\tau_{\mathbf1}\wedge\tau_{\mathbf0}<\infty
)=\P^w(\tau_{\mathbf1}<\infty)+\P^w(\tau_{\mathbf0}<\infty),
\]
where we write $\tau_{\eta}$ for the first hitting time of $\eta$.
Observe that $\P^w(\tau_{\mathbf1}<\infty)$ is independent of the
auxiliary number $\lambda>0$. Indeed, the holding time of each
configuration $\eta\neq\mathbf1,\mathbf0$ is finite and the
probability of transition from $\eta$ to $\eta^x$ at the end of the
holding time is given by
\[
\frac{c^w(x,\eta)}{\sum_{y\in\mathsf V} c^w(y,\eta)},
\]
which is independent of $\lambda>0$.

We can estimate the equilibrium probability $\P^w(\tau_{\mathbf
1}<\infty)$ by a ``harmonic sampler'' of the voter model from finite time.
Let $p_1(\xi)$ be the weighted average of $1$'s in the vertex set
%
%e27 #&#
\begin{equation}
p_1(\eta)=\sum_x \eta(x)\pi(x),\label{eqp1}
\end{equation}
where $\pi(x)$ is the invariant distribution of the (nearest-neighbor)
random walk on $G$ and is given by (\ref{eqinv}).
Since $p_1(\mathbf1)=1-p_1(\mathbf0)=1$ and the chain is eventually
trapped at $\mathbf1$ or $\mathbf0$, it follows from dominated
convergence that
%
%e28 #&#
\begin{equation}
\lim_{n\to\infty}\Bbb E^w[p_1(\xi_n)]=\Bbb P^w(\tau_{\mathbf
1}<\infty).\label{eqlp}
\end{equation}
On the other hand, the function $p_1$ is harmonic for the voter model
\begin{eqnarray*}
\Bbb E_\eta[p_1(\xi_1)]%&\frac{\lambda}{N}\sum_{\xi(x)=1}(p_1(
%&+
%&+p_1(\xi)(1-\frac{\lambda}{N}\sum_{\xi(x)=1}\frac{\#\{y\sim%x;
&=&p_1(\eta)+\frac{\lambda}{N\cdot2\# \mathsf E}\\
&&\hspace*{36pt}{}\times\biggl(\sum_{\eta(x)=0}\!\#\{y\sim x;\eta(y)=1\}-\!\sum
_{\eta(x)=1}\!\#\{y\sim x;\eta(y)=0\}\!\biggr)\\
&=&p_1(\eta).
\end{eqnarray*}
In particular, (\ref{eqlp}) applied to $w=0$ entails
%
%e29 #&#
\begin{equation}
\P_\eta(\tau_{\mathbf1}<\infty)=p_1(\eta)=\frac{\sum_{\eta
(x)=1}d (x)}{2\cdot\#\mathsf E},\label{eqpiREG}
\end{equation}
where the last equality follows from the explicit form (\ref{eqinv})
of $\pi$.\vadjust{\goodbreak}

%re4 #&#
\begin{rmk}\label{rmkharm}
Since every harmonic function $f$ for the voter model satisfies
\[
f(\eta)\equiv\Bbb E_\eta[f(\xi_{\tau_{\mathbf1}\wedge\tau
_{\mathbf0}})],
\]
(\ref{eqpiREG}) implies that the vector space of harmonic functions
is explicitly characterized as the span of the constant function $1$
and $p_1$. Recall also the foregoing display gives a construction of
any harmonic function with preassigned values at $\mathbf1$ and
$\mathbf0$. (See, e.g., Chapter 2 in~\cite{AFIPS}.) %\qed
\end{rmk}

%s3 #&#
\section{Expansion}\label{secexp}
We continue to study the discrete-time voter model perturbations
defined in Section~\ref{secbs}.
For each $w\in[0,w_0]$, consider the signed kernel
\[
K^w=P^w-P,
\]
which measures the magnitude of perturbations of transition probabilities.
We also define a nonnegative kernel $|K^w|$ by $|K^w|(\eta,\widetilde
{\eta})=|K^w(\eta,\widetilde{\eta})|$.

%le5 #&#
\begin{lem}\label{lemKwest}
For any $w\in[0,w_0]$ and any $f\dvtx \{1,0\}^{\mathsf V}\lra\Bbb R$, we have
%
%e30 #&#
\begin{equation}
K^wf(\eta)=\frac{\lambda}{N}\sum_x\bigl[wh_{1-\eta(x)}(x,\eta
)+w^2g_w(x,\eta)\bigr][
f(\eta^x)-f(\eta)]\label{eqKw}
\end{equation}
and
%
%e31 #&#
\begin{equation}
\||K^w|f\|_\infty\leq4C_0w\|f\|_\infty,\label{ineqKw}
\end{equation}
where $C_0$ is the constant in \textup{(A1)}.
\end{lem}
\begin{pf} %{Proof of Lemma~\ref{lemKwest}}
We notice that for any $\eta$ and any $x$,
\begin{eqnarray*}
K^w(\eta,\eta^x)&=&\frac{\lambda}{N}\bigl[wh_{1-\eta(x)}(x,\eta
)+w^2g_w(x,\eta)\bigr],\\
K^w(\eta,\eta)&=&-\sum_{x}\frac{\lambda}{N}\bigl[wh_{1-\eta
(x)}(x,\eta)+w^2g_w(x,\eta)\bigr],
\end{eqnarray*}
by the definitions of $c^w$ and $P^w$. Our assertions (\ref{eqKw})
and (\ref{ineqKw}) then follow at one stroke.
\end{pf}

Using the signed kernel $K^w$, we can rewrite every $T$-step transition
probability $P^w_T$ of the voter model perturbation as
%
%e32 #&#
%e33 #&#
%e34 #&#
\begin{eqnarray}\label{eqPwT}
P^w_T(\eta_0,\eta_T)&=&\sum_{\eta_1,\ldots,\eta_{T-1}}P^w(\eta
_0,\eta_1)
P^w(\eta_1,\eta_2)\cdots P^w(\eta_{T-1},\eta_T)\nonumber\\
&=&\sum_{\eta_1,\ldots,\eta_{T-1}}(P+K^w)(\eta_0,\eta_1)\cdots
(P+K^w)(\eta_{T-1},\eta_T)\\
&=&P_T(\eta_0,\eta_T)+\sum_{n=1}^T\sum_{\mathbf j\in\mathcal
I_T(n)}\sum_{\eta_1,\ldots,\eta_{T-1}}\Delta^{w,\mathbf j}_T(\eta
_0,\ldots,\eta_T).\nonumber
\end{eqnarray}
Here, $\mathcal I_T(n)$ is the set of strictly increasing $n$-tuples
with entries in $\{1,\ldots,T\}$, and for $\mathbf j=(j_1,\ldots
,j_n)\in\mathcal I_T(n)$
%
%e35 #&#
\begin{equation}
\Delta^{w,\mathbf j}_T(\eta_0, \ldots,\eta_T)\label{eqDelta}
\end{equation}
is the signed measure of the path $(\eta_0, \eta_1, \ldots, \eta
_T)$ such that the transition from $\eta_{r}$ to $\eta_{r+1}$ is
determined by $K^w(\eta_r,\eta_{r+1})$ if $r+1$ is one of the
(integer-valued) indices in $\mathbf j$ and is determined by $P(\eta
_r,\eta_{r+1})$ otherwise. For convenience, we set for each $\mathbf
j\in\mathcal I_T(n)$
\[
Q_T^{w,\mathbf j}(\eta_0,\eta_T)=\sum_{\eta_1,\ldots,\eta
_{T-1}}\Delta^{w,\mathbf j}_T(\eta_0,\ldots,\eta_T)
\]
as the $T$-step transition signed kernel, and we say $Q^{w,\mathbf
j}_T$ has $n$ \textit{faults} (up to time $T$) and $\mathbf j$ is its \textit{fault sequence}. Then by (\ref{eqPwT}), we can write for any $f\dvtx \{
1,0\}^{\mathsf V}\lra\Bbb R$,
%
%e36 #&#
%e37 #&#
\begin{eqnarray}\label{eqexp}
\Bbb E^w_{\eta_0}[f(\xi_T)]&=&\Bbb E_{\eta_0}[f(\xi_T)]+\sum
_{n=1}^T\sum_{\mathbf j\in\mathcal I_T(n)}\sum_{\eta_1,\ldots,\eta
_T}\Delta^{w,\mathbf j}_T(\eta_0,\ldots,\eta_T)f(\eta_T)
\nonumber
\\[-8pt]
\\[-8pt]
\nonumber
&=&\Bbb E_{\eta_0}[f(\xi_T)]+\sum_{n=1}^T\sum_{\mathbf j\in
\mathcal I_T(n)}Q^{w,\mathbf j}_Tf(\eta_0).
\end{eqnarray}

Write $\mathcal I(n)\equiv\mc I_\infty(n)$ for the set of strictly
increasing $n$-tuples with entries in~$\Bbb N$. We now state the key
result which in particular offers an expansion of fixation probabilities.

%pr6 #&#
\begin{prop}\label{propexp}
Recall the parameter $w_0>0$ in the definition of the voter model perturbations.
There exists $w_1\in(0,w_0]$ such that for any harmonic function~$f$
for the voter model,
%
%e38 #&#
\begin{equation}
\quad f(\mathbf1)\P^w_{\eta}(\tau_{\mathbf1}<\infty)+f(\mathbf0)\P
_{\eta}^w(\tau_{\bsm0}<\infty) =f(\eta)+\sum_{n=1}^\infty\sum
_{\mathbf j\in\mathcal I(n)}Q^{w,\mathbf j}_{j_n}f (\eta),\label{eqseries}
\end{equation}
where the series converges absolutely and uniformly in $w\in[0,w_1]$
and in $\eta\in\{1,0\}^{\mathsf V}$.
\end{prop}

%re7 #&#
\begin{rmk}
(i) There are alternative perspectives to state the conclusion of
Proposition~\ref{propexp}. Thanks to Remark~\ref{rmkharm}
and the fact $Q^{w,\mathbf j}_T\1\equiv0$, it is equivalent to the
validity of the same expansion for $p_1$ defined in (\ref{eqp1})
(for any small $w$). By Remark~\ref{rmkharm} again, it is also
equivalent to an analogous expansion of any linear combination of the
two fixation probabilities under~$\P^w$.\vspace*{-6pt}
\begin{longlist}[(ii)]
\item[(ii)] The series expansion (\ref{eqseries}) has the flavor of a Taylor
series expansion in $w$, as hinted by Lemma~\ref{lemest}. \vadjust{\goodbreak}%\qed
\end{longlist}
\end{rmk}

The proof of Proposition~\ref{propexp} is obtained by passing $T$ to
infinity for both sides of~(\ref{eqexp}). This immediately gives the
left-hand side of (\ref{eqseries}) thanks to our assumption (\ref
{A4}). There are, however, two technical issues when we handle the
right-hand sides of (\ref{eqexp}). The first one is minor and is the
dependence on $T$ of the summands $Q^{w,\mathbf j}_Tf(\eta_0)$. For
this, the harmonicity of $f$ implies that such dependence does not
exist, as asserted in Lemma~\ref{lemmg}. As a result, the remaining
problem is the absolute convergence of the series on the right-hand
side of (\ref{eqseries}) for any small parameter $w>0$. This is
resolved by a series of estimates in Lemmas~\ref{lemexp}, \ref
{lemest} and finally Lemma~\ref{lemseries}.

%le8 #&#
\begin{lem}\label{lemmg}
For any harmonic function $f$ for the voter model, any $T\geq1$, and
any $\mathbf j\in\mc I_T(n)$,
\[
Q^{w,\mathbf j}_Tf(\eta_0)\equiv Q^{w,\mathbf j}_{j_n}f(\eta_0),
\]
where we identify $\mathbf j\in\mc I_{j_n}(n)$ in the natural way.
\end{lem}
\begin{pf}
This follows immediately from the martingale property of a harmonic
function $f$ for the voter model and the definition of the signed
measures $\Delta^{w,\mathbf j}_{T}$ in~(\ref{eqDelta}).
\end{pf}
%

%le9 #&#
\begin{lem}\label{lemexp}
There exist $C_1=C_1(G)\geq1$ and $\delta=\delta(G)\in(0,1)$ such that
%
%e39 #&#
\begin{equation}
\sup_{\eta\neq\mathbf1,\mathbf0}\Bbb P_\eta(\xi_n\neq\mathbf
1,\mathbf0)\leq C_1\delta^n\qquad \mbox{for any $n\geq1$}.\label{ineqexp}
\end{equation}
\end{lem}
\begin{pf}
Recall that the voter model starting at any arbitrary state is
eventually trapped at either $\mathbf1$ or $\mathbf0$.
By identifying $\mathbf1$ and $\mathbf0$, we deduce (\ref{ineqexp})
from some standard results of nonnegative matrices, for suitable
constants $C_1=C_1(G)\geq1$ and $\delta\in(0,1)$. (See, e.g., \cite
{AAPQ}, % Section I.6 and
Lemma I.6.1 and Proposition~I.6.3.)
\end{pf}

%le10 #&#
\begin{lem}\label{lemest}
Let $C_1=C_1(G)$ and $\delta=\delta(G)$ be the constants in
Lemma~$\ref{lemexp}$, and set $C=C(G,C_0)=\max(4C_0,C_1)$. Then
for any $\mathbf j\in\mathcal I(n)$, any $w\in[0,w_0]$ and any
harmonic function $f$ for the voter model,
%
%e40 #&#
\begin{equation}
\|Q^{w,\mathbf j}_{j_n}f\|_\infty\leq\|f\|_\infty
w^nC^{2n}\delta^{j_n-n}.\label{eqest}
\end{equation}
%
%with the convention $Q^{w,\varnothing}_jf=Q^{w,j}_jf$ for $j\in\mc
%I(1)$.
\end{lem}
\begin{pf}
Without loss of generality, we may assume $\|f\|_\infty=1$. By definition,
%
%e41 #&#
\begin{equation}\label{eqaux}
Q^{w,\mathbf j}_{j_n}f(\eta_0)=\sum_{\eta_1,\ldots,\eta
_{j_n}}\Delta^{w,\mathbf j}_{j_n}(\eta_0,\ldots,\eta_{j_n})f(\eta_{j_n}).
\end{equation}
If
$\Delta^{w,\mathbf j}_{j_n}(\eta_0,\ldots,\eta_{j_n})$ is nonzero, then
none of the $\eta_0,\ldots,\eta_{j_n-1}$ is $\mathbf1$ or $\mathbf
0$. Indeed, if some of the $\eta_i$, $0\leq i\leq j_n-1$, is $\mathbf
1$, then $\eta_{i+1},\ldots,\eta_{j_n}$ can only be $\mathbf1$ by~(\ref{A3}) and therefore $K^wf(\eta_{j_n-1})=0$. This is a contradiction.
Similarly, we cannot have $\eta_i=\mathbf0$ for some $0\leq i\leq
j_n-1$. Hence, the nonvanishing summands of the right-hand side of~(\ref{eqaux}) range over $\Delta^{w,\mathbf j}_{j_n}(\eta
_1,\ldots,\eta_{j_n})f(\eta_{j_n})$ for which none of the
$\eta_{j_1},\ldots,\eta_{j_n-1}$ is $\mathbf1$ or $\mathbf0$. With
$\eta_0$ fixed, write $\Delta^{w,\mathbf j'}_{U,\eta_0}$ for the
signed measure $\Delta^{w,\mathbf j'}_U$ restricted to paths starting
at $\eta_0$.
Thus we get from (\ref{eqaux}) that
\[
Q^{w,\mathbf j}_{j_n}f(\eta_0)=\Delta^{w,\mathbf j}_{j_n,\eta
_0}\bigl[f(\xi_{j_n})\1_{[\xi_1,\ldots,\xi_{j_n-1}\neq
\mathbf1,\mathbf0]}\bigr].
\]
Here, our usage of the compact notation on the right-hand side is
analogous to the convention in the modern theory of stochastic processes.
Recall that\vspace*{-1pt} $|K^w|$ stands for the kernel $|K^w|(\eta,\tilde{\eta
})=|K^w(\eta,\tilde{\eta})|$, and put $|\Delta|^{w,\mathbf
j'}_{\eta_0,U}$ for the measure on paths $(\eta_0,\ldots,\eta_U)$
obtained by replacing all the $K^w$ in $\Delta^{w,\mathbf j'}_{U,\eta
_0}$ by $|K^w|$. Since $\|f\|_\infty=1$, the foregoing display implies
\begin{eqnarray*}
|Q^{w,\mathbf j}_{j_n}f(\eta_0)|&\leq&|\Delta|^{w,\mathbf
j}_{j_n,\eta_0}(\xi_1,\ldots,\xi_{j_n-1}\neq\mathbf1,\mathbf0)\\
&\leq&|\Delta|^{w,(j_1,\ldots,j_{n-1})}_{j_{n-1},\eta_0}(\xi
_1,\ldots,
\xi_{j_{n-1}-1}
\neq\mathbf1,\mathbf0)\\
&&{}\times\Bigl(\sup_{\eta\neq\mathbf1,\mathbf0}\Bbb P_{\eta
}(\xi_{j_n-j_{n-1}-1}\neq\mathbf1,\mathbf0)\Bigr) \||K^w|\1
\|_\infty
\end{eqnarray*}
with $j_0=0$,
where $\sup_{\eta\neq\mathbf1,\mathbf0}\Bbb P_{\eta}(\xi
_{j_n-j_{n-1}-1}\neq\mathbf1,\mathbf0)$ bounds the measure of the
yet ``active'' paths from $j_{n-1}$ to $j_n-1$ and $\||K^w|\1\|
_\infty$ bounds the measure of the transition from $j_n-1$ to $j_n$.
Iterating the last inequality, we get
%
%e42 #&#
\begin{eqnarray}\label{eqaux1}
|Q^{w,\mathbf j}_{j_n}f(\eta_0)|&\leq&\||K^w|\1\|^n_\infty
\prod_{r=1}^n\Bigl(\sup_{\eta\neq\mathbf1,\mathbf0}\Bbb P_{\eta}
(\xi_{j_r-j_{r-1}-1}\neq\mathbf1,\mathbf0)\Bigr)
\nonumber
\\[-8pt]
\\[-8pt]
\nonumber
&\leq&(4C_0)^nw^n\prod_{r=1}^n\Bigl(\sup_{\eta\neq\mathbf
1,\mathbf0}\Bbb P_{\eta}(\xi_{j_r-j_{r-1}-1}\neq\mathbf1,\mathbf
0)\Bigr),
\end{eqnarray}
where the last inequality follows from Lemma~\ref{lemKwest}.
%In addition, by (\ref{A1}), we can modify the constant $C$ such that
%(cf. the proof of Lemma~\ref{lemKwest}).
Since $\sum_{r=1}^n(j_r-j_{r-1}-1)=j_n-n$ and $C=\max(4C_0,C_1)$,
Lemma~\ref{lemexp} applied to the right-hand side of (\ref{eqaux1}) gives
%
%e43 #&#
\begin{equation}
|Q^{w,\mathbf j}_{j_n}f(\eta_0)|\leq w^n(C^2)^{n}\delta
^{j_n-n}.\label{ineqcoreest}
\end{equation}
The proof of (\ref{eqest}) is complete.
\end{pf}

%le11 #&#
\begin{lem}\label{lemseries}
Recall the constants $C=C(G,C_0)$ and $\delta=\delta(G)$ in
Lemmas~\ref{lemest} and~\ref{lemexp}, respectively. There
exists $w_1\in(0,w_0]$ such that
%
%e44 #&#
\begin{equation}
\sum_{n=1}^\infty\sum_{\mathbf j\in\mc I(n)}w_1^n(C^2)^{n}\delta
^{j_n-n}<\infty.\label{eqseriesconv}
\end{equation}
\end{lem}
\begin{pf}
Observe that every index in $\bigcup_{n=1}^\infty\mathcal I(n)$ can
be identified uniquely by the time of the last fault and the fault
sequence before the time of the last fault. Hence, letting $S$ denote
the time of the last fault and $m$ the number of faults within $\{
1,\ldots,S-1\}$, we can write for any $w>0$
%
%e45 #&#
\begin{equation}
\qquad\sum_{n=1}^\infty\sum_{\mathbf j\in\mc I(n)}w^n(C^2)^n\delta^{j_n-n}=
\sum_{S=1}^\infty\sum_{m=0}^{S-1}\pmatrix{S-1\vspace*{2pt}\cr
m}w^{m+1}(C^2)^{m+1}\delta^{S-m-1}.\label{eqaux01}
\end{equation}
For each $S$, write
%
%e46 #&#
\begin{equation}
\qquad\sum_{m=0}^{S-1}\pmatrix{S-1\cr m}w^{m+1}(C^2)^{m+1}\delta
^{S-m-1}=wC^2\delta^{S-1}(1+wC^2\delta^{-1})
^{S-1}.\label{eqaux7}
\end{equation}
With $\delta\in(0,1)$ fixed, we can choose $w_1\in(0,w_0]$ small
such that
\[
C^2(1+w_1C^2\delta^{-1})^{S-1}\leq\biggl(\frac{1}{\sqrt{\delta
}}\biggr)^{S-1}
\]
for any large $S$. Apply the foregoing inequality for large $S$ to the
right-hand side of (\ref{eqaux7}) with $w$ replaced by $w_1$. This gives
\[
%w_1C^2\delta^{S-1}(1+wC^2\delta^{-1})^{S-1}
\sum_{m=0}^{S-1}\pmatrix{S-1\cr m}w_1^{m+1}(C^2)^{m+1}\delta^{S-m-1}\leq
w_1\bigl(\sqrt{\delta}\bigr)^{S-1},
\]
where the right-hand side converges exponentially fast to $0$ as
$\delta<1$. By (\ref{eqaux01}), the asserted convergence of the
series in (\ref{eqseriesconv}) now follows.
\end{pf}

\begin{pf*}{Proof of Proposition~\ref{propexp}}
We pick $w_1$ as in the statement of Lem\-ma~\ref{lemseries}.
By Lemma~\ref{lemest} and the choice of $w_1$, the series in (\ref
{eqseries}) converges absolutely and uniformly in $w\in[0,w_1]$ and
$\eta\in\{1,0\}^{\mathsf V}$. %; for the other indices $S$, we still
%have
%wC^2\delta^{S-1}((1+wC^2\delta^{-1})^{S-1}-1)\leq w^2C_1
%for some constant $C_1$ by the mean value theorem. Hence, (

By (\ref{eqexp}) and dominated convergence, it remains to show that
\[
\lim_{T\to\infty}\sum_{n=1}^T\sum_{\mathbf j\in\mc
I_T(n)}Q^{w,\mathbf j}_Tf(\eta_0)= \sum_{n=1}^\infty\sum_{\mathbf
j\in\mc I(n)}Q^{w,\mathbf j}_{j_n}f(\eta_0).
\]
To see this, note that by Lemma~\ref{lemmg}, we can write
\[
\sum_{n=1}^T\sum_{\mathbf j\in\mathcal I_T(n)}Q^{w,\mathbf
j}_Tf(\eta_0)=\sum_{n=1}^T\mathop{\sum_{\mathbf j\in
\mc I(n)}}_{j_n\leq T} Q^{w,\mathbf j}_{j_n}f(\eta_0),
\]
where the right-hand side is a partial sum of the infinite series in
(\ref{eqseries}). The validity of (\ref{eqseries}) now follows from
the absolute convergence of the series in the same display. The proof
is complete.
\end{pf*}

%The expansion (\ref{eqseries}) gives finite order approximations in
%$w$ of the fixation probability $\P^w(\tau_{\mathbf1}<\infty)$.
%Indeed, taking $\gamma_{\bs1}=1-\gamma_{\bs0}=1$ and using the bound
%in (\ref{eqest}), we see each partial sum
%has order $w^n$ for each $n$. This suggests we can estimate $\P^w(
%that we are able to calculate these partial sums. In the sequel, we
%only consider the first order approximation, namely $T_0=1$, of $\P^w(

For the convenience of subsequent applications, we
consider from now on the continuous-time Markov chain $(\xi_s)$ with
rates given by (\ref{A2}). We can define this chain $(\xi_s)$ from
the discrete-time Markov chain $(\xi_n)$ by
\[
\xi_s=\xi_{M_s},
\]
where $(M_s)$ is an independent Poisson process with $\Bbb E[M_s]=\frac
{sN}{\lambda}$. (Recall our time scale convention: ``$n$'' for the
discrete time scale and ``$s$'' for the continuous time scale.)
Under this setup, the potential measure of $(\xi_s)$ and the potential
measure of $(\xi_n)$ are linked by
%
%e47 #&#
\begin{equation}
\Bbb E\biggl[\int_0^\infty f(\xi_s)\,ds\biggr]=\frac{\lambda}{N}\Bbb
E\Biggl[\sum_{n=0}^\infty f(\xi_n)\Biggr]\label{eqsumint}
\end{equation}
for any nonnegative $f$.
In addition, the fixation probability to $\mathbf1$ for this
continuous-time Markov chain $(\xi_s)$ is the same as that for the
discrete-time chain $(\xi_n)$. (See the discussion after
Proposition~\ref{propiu}.)

We now state a first-order approximation of $\P^w(\tau_{\bsm1}<\infty
)$ by the voter model. Recall the difference kernel $D$ defined by
(\ref{eqD}) with $h_i$ as in (\ref{A1}) and the $\pi$-expectation
$\overline{D}$ defined by (\ref{eqHbar}).

%th3.1 #&#
\begin{thmm}\label{thmTaylor}
Let $\nu$ be an arbitrary distribution on $\{1,0\}^{\mathsf V}$. Then
as $w\lra0+$, we have
%
%e48 #&#
\begin{equation}
\Bbb P^w_\nu(\tau_{\mathbf1}<\infty)=\Bbb P_\nu(\tau_{\mathbf
1}<\infty)+w\int_0^\infty\Bbb E_\nu[\overline{D}(\xi
_s)]\,ds+O(w^2).\label{ineqcrt}
\end{equation}
Here, the convention for the function $O(w^2)$ is as in Theorem~\ref
{thmi1}. Moreover, $\overline{D}(\xi_s)\in L_1(d\P_\nu\otimes ds)$.
\end{thmm}
\begin{pf}
It suffices to prove the theorem for $\nu=\delta_{\eta_0}$ for any
$\eta_0\in\{1,0\}^{\mathsf V}$. Recall that the function
$p_1$ defined by (\ref{eqp1}) is harmonic for the voter model, and
hence, the expansion (\ref{eqseries}) applies. By (\ref{eqest}) and
Lemma~\ref{lemseries}, it is plain that
%
%e49 #&#
\begin{equation}
\sum_{n=2}^\infty\sum_{\mathbf j\in\mc I(n)}Q^{w,\mathbf
j}_{j_n}p_1(\eta_0)=O(w^2).\label{eqaux3}
\end{equation}
We identify each $\mathbf j\in\mc I(1)$ as $\mathbf j=(j)=j$ and look
at the summands $Q^{w,j}_jp_1$.
Write $\Bbb E^{w,j}$ for the expectation of the time-inhomogeneous
Markov chain where the transition of each step is governed by $P$
except that the transition from $j-1$ to $j$ is governed by $P^w$. Then
%
%e50 #&#
\begin{eqnarray}\label{eqaux5}
Q^{w,j}_jp_1(\eta_0)&=&\Bbb E^{w,j}_{\eta_0}[p_1(\xi_j)]-\Bbb
E_{\eta_0}[p_1(\xi_j)]\nonumber\\
&=&\Bbb E_{\eta_0}\bigl[ \Bbb E^w_{\xi_{j-1}}[p_1(\xi_1)]-\Bbb
E_{\xi_{j-1}}[p_1(\xi_1)]\bigr]
\nonumber
\\[-8pt]
\\[-8pt]
\nonumber
&=&\Bbb E_{\eta_0}[K^wp_1(\xi_{j-1})]\\
&=&\Bbb E_{\eta_0}[K^wp_1(\xi_{j-1});\tau_{\mathbf1}\wedge\tau
_{\mathbf0}\geq j],\nonumber
\end{eqnarray}
where the last equality follows from the definition of $K^w$ and the
fact that $\mathbf1$ and $\mathbf0$ are both absorbing.
Moreover, we deduce from Lemma~\ref{lemKwest} that
%
%e51 #&#
\begin{equation}
K^wp_1(\eta)=\frac{\lambda}{N}w\overline{D}(\eta)+\frac{\lambda}{N}w^2
\overline{G}_w(\eta),\label{eqaux6}
\end{equation}
where
\[
G_w(x,\eta)=g_w(x,\eta)\bigl(1-2\eta(x)\bigr).
\]
Note that $\Bbb E_{\eta_0}[\tau_{\mathbf1}\wedge\tau_{\mathbf
0}]<\infty$ by Lemma~\ref{lemexp}. Hence, by (\ref{eqaux5}) and
(\ref{eqaux6}), we deduce that
%
%e52 #&#
\begin{eqnarray}\label{eqaux4}
\sum_{j=1}^\infty Q^{w,j}_jp_1(\eta_0)&=&\frac{\lambda w}{N}\sum
_{j=1}^\infty\Bbb E_{\eta_0}[\overline{D}(\xi_{j-1})
]+O(w^2)
\nonumber
\\[-8pt]
\\[-8pt]
\nonumber
&=&w\Bbb E_{\eta_0}\biggl[\int_0^\infty\overline{D}(\xi_s)\,ds
\biggr]+O(w^2),
\end{eqnarray}
where the last equality follows from (\ref{eqsumint}). Moreover,
$\overline{D}(\xi_s)\in L_1(d\P_{\eta_0}\otimes ds)$.
The approximation (\ref{ineqcrt}) for each $w\leq w_2$ for some small
$w_2\in(0,w_1]$ now follows from (\ref{eqaux3}) and (\ref{eqaux4})
applied to the expansion (\ref{eqseries}) for $p_1$. The proof is complete.
\end{pf}

\section{First-order approximations}\label{secbc}
In this section, we give the proof of Theorem~\ref{thmi1}. We consider only
regular graphs throughout this section. (Recall that a graph is regular
if all vertices have the same number of neighbors.)

As a preliminary, let us introduce the convenient notion of Bernoulli
transforms and discuss its properties.
For each $u\in[0,1]$, let $\mu_u$ be the Bernoulli product measure on
$\{1,0\}^{\mathsf V}$ with density $\mu_u(\xi(x)=1)=u$. For any
function $f\dvtx \{1,0\}^{\mathsf V}\lra\Bbb R$, define the \textit{Bernoulli
transform} of $f$ by
%
%e53 #&#
\begin{eqnarray}\label{eqBf}
\B f(u):=\int f\,d\mu_u=\sum_{n=0}^N\biggl[\sum_{\eta: \#\{x;\eta
(x)=1\}=n}
f(\eta)\biggr]u^n(1-u)^{N-n},
\nonumber
\\[-8pt]
\\[-8pt]
\eqntext{u\in[0,1].}
\end{eqnarray}
The Bernoulli transform of $f$ uniquely determines the coefficients
\[
A_f(n):= \sum_{\eta: \#\{x;\eta(x)=1\}=n}f(\eta),\qquad 0\leq n\leq N.
\]
Indeed, $\B f(0)=f(\mathbf0)=A_f(0)$ and for each $1\leq n\leq N$,
\[
A_f(n)=\lim_{u\da0+}\frac{1}{u^n}\Biggl(\B I(u)-\sum
_{i=0}^{n-1}u^i(1-u)^{N-i}A_f(i)\Biggr).
\]

The Bernoulli transform $\B f(u)$ is a polynomial $\sum
_{i=0}^{N}\alpha_i u^i$ of order at most~$N$. Let us
invert the coefficients $A_f(n)$ from $\alpha_i$ by basic
combinatorics. By the binomial theorem,
\[
u^i=\sum_{n=0}^{N-i}\pmatrix{N-i\cr n}u^{i+n}(1-u)^{N-i-n}.
\]
Hence, summing over $i+n$, we have
\[
\sum_{i=0}^N\alpha_iu^i=\sum_{n=0}^{N}\Biggl[\sum_{i=0}^n\alpha
_i\pmatrix{N-i\cr n-i}\Biggr] u^n(1-u)^{N-n},
\]
and the uniqueness of the coefficients $A_f$ implies
%
%e54 #&#
\begin{equation}
A_f(n)=\sum_{i=0}^n\alpha_i\pmatrix{N-i\cr n-i},\qquad 0\leq n\leq
N.\label{eqAf}
\end{equation}
As a corollary, we obtain
%
%e55 #&#
\begin{equation}
\int f\,d\mathbf u_n=\frac{1}{{N\choose n}}A_f(n)=\frac{1}{{N\choose
n}}\sum_{i=0}^n\alpha_i\pmatrix{N-i\cr n-i},\qquad 1\leq n\leq N-1,\label{eqform}
\end{equation}
if we regard $\mathbf u_n$, the uniform distribution on the set of
configurations with precisely~$n$ many $1$'s, as a measure on $\{1,0\}
^{\mathsf V}$ in the natural way.

%We note that
%&+\ldots+u^{N-1}(1-u)\sum_{|\{x;;\eta(x)=1\}|=N-1}f(\eta)+u^Nf(
%This implies %we can invert the function
%n\lmt\sum_{|\{x;\eta(x)=0\}|=n}f(\eta),
%from the Bernoulli transform $\B f$. Indeed,
%$f(\mathbf0)=\B f(0)$ and
%%f(0)}{u}.\label{eqB1}
We will specialize the application of Bernoulli transforms to the function
\[
I(\eta):=\int_0^\infty\Bbb E_\eta[\overline{D}(\xi_s)]\,ds.
\]
To obtain the explicit approximations (up to the first order) asserted
in Theorem~\ref{thmi1}, we need to compute by Theorem~\ref{thmTaylor} the
$0$-potentials $\int_0^\infty\E_{\mathbf u_n}[\overline
{D}(\xi_s)]\,ds$ for $1\leq n\leq N-1$ under each updating
mechanism. On the other hand, we will see that the Bernoulli transform
of each $0$-potential $I$ is analytically tractable and
%
%e56 #&#
\begin{equation}
\B I(u)=\Gamma u(1-u)\label{eqBfSpec}
\end{equation}
for some explicit constant $\Gamma$. Note that we have
$A_I(N)=A_I(0)=0$ for the updating mechanisms under consideration.
Hence, the formula (\ref{eqform}) entails
%
%e57 #&#
\begin{equation}
\int_0^\infty\Bbb E_{\mathbf u_n}[\overline{D}(\xi_s)
]\,ds=\frac{\Gamma n(N-n)}{N(N-1)},\qquad 1\leq n\leq N-1,\label{eqBIu}
\end{equation}
since ${N-1\choose n-1}-{N-2\choose n-2}={N-2\choose n-1}$ for $n\geq
2$ and ${N-1\choose0}={N-2\choose0}=1$.
%By Corollary~\ref{corselection}, whether selection strongly favors or
%strongly opposes cooperation depends only on the sign of the
%derivative of $\B I$ at $u=0$, if the location of the cooperative
%mutant is uniformly distributed on the vertices.
%Variants will also be discussed in Example~\ref{egKn} and Example~
%Precisely, under the birth--death updates, we shall show that whenever
%a payoff matrix (\ref{eqPi}) and a degree $k$ are both given, then
%there exists a positive integer $N_0$ such that for all $k$-regular
%symmetric graphs, selection favors $1$ replacing $0$ whenever $b-

%s4.1 #&#
\subsection{\texorpdfstring{Proof of Theorem \protect\ref{thmi1}(1)}{Proof of Theorem 1(1)}}\label{subsecdbuT}
Assume that the graph $G$ is $k$-regular.
Recall that the death--birth updating defines a Markov chain of voter
model perturbation satisfying (\ref{A1})--(\ref{A4}) by
Proposition~\ref{propdbu} under weak selection. The functions $h_i$
in~(\ref{A1}) for this updating are given by (\ref{eqdbh}). Hence,
the difference kernel $D$ is given by $D(x,\xi)=h_1(x,\xi)$, and for
any $x\in\mathsf V$,
%
%e58 #&#
%e59 #&#
\begin{eqnarray}\label{eqaux2}
\frac{1}{k}\Bbb E_{\mu_u}[D(x,\xi_s)]
%=&\frac{1}{k}\Bbb E_{\mu_u}[h_1(o,\xi_s)]\nonumber\\
&=&-(b+c)\Bbb E_{\mu_u}[f_0f_1(x,\xi_s)]+b\Bbb E_{\mu
_u}[f_{00}(x,\xi_s)]\nonumber\\
&&{}+b\Bbb E_{\mu_u}[f_0f_{11}(x,\xi_s)]-b\Bbb E_{\mu
_u}[f_0f_{00}(x,\xi_s)]
\nonumber
\\[-8pt]
\\[-8pt]
\nonumber
&=&-(b+c)\Bbb E_{\mu_u}[f_0f_1(x,\xi_s)]
+b\Bbb E_{\mu_u}[f_0f_{11}(x,\xi_s)]\\
&&{}+b\Bbb E_{\mu_u}[f_1f_{00}(x,\xi_s)].\nonumber
\end{eqnarray}

In analogy to the computations in~\cite{CDPVMP} for coalescing
probabilities, we resort to the duality between voter models and
coalescing random walks for the right-hand side of (\ref{eqaux2}).
Let $\{B^x;x\in\mathsf V\}$ be the rate-$1$ coalescing random walks on
$G$, where~$B^x$ starts at $x$. The random walks move independently of
each other until they meet another and move together afterward. The
duality between the voter model and the coalescing random walks is
given by
\[
\Bbb P_\eta\bigl(\xi_s(x)=i_x,x\in Q\bigr)=\Bbb P\bigl(\eta
(B^x_s)=i_x, x\in Q\bigr)
\]
for any $Q\subset\mathsf V$ and $(i_x;x\in Q)\in\{1,0\}^Q$. (See
Chapter V in~\cite{LIPS}.)
Introduce two independent discrete-time random walks $(X_n;n\geq0)$
and $(Y_n;n\geq0)$ starting at the same vertex, both independent of $\{
B^x;x\in\mathsf V\}$.
%e_1\stackrel{(\rm d)}{=}X^o_{J_1^o},
%(e_2,e_3)\stackrel{(\rm d)}{=}(X^o_{J_1^o},X^o_{J_2^o}),
%and we assume $e_1,e_2,e_3$ are independent of $X^x$, $x\in\mathsf V$.
%Then we have
%%dx,Y_1\in dy)\Bbb P_\eta(\xi_s(x)=0,\xi_s(y)=1)\\
%=&\int\mu_u(d\eta)\int\Bbb P(X_1\in dx, Y_1\in dy)\Bbb%P(
%=&\int\mu_u(d\eta)\Bbb%P(\eta(B^{X_1}_s)=0, \eta(B^{Y_1}_s)=1
%),
%%=&\int\Bbb P(B^{e_1}_s\in dx,B^{e_2}_s\in dy)\1_{x\neq y}u(1-u)\\
%%=&u(1-u)p_s(e_1|e_2).
%and similar equalities hold for $\Bbb E_{\mu_u}[f_0f_{11}(x,\xi_s)]$
%and $\Bbb E_{\mu_u}[f_1f_{00}(x,\xi_s)]$.
%Similar computations apply to other terms on the right-hand side of (
%We now apply the convenient notation in~\cite{CDPVMP}. %Set
%(\eta(B^{e_1}_s)=0, \eta(B^{e_2}_s)=1),
%where $\dxi=1-\xi$ and
%Define $M_{x,y}=\inf\{t\in\Bbb R_+;Z^x_t=Z^y_t\}$ for the first
%meeting time of the independent walks $Z^x$ and $Z^y$.
Fix $x$ and assume that the chains $(X_n)$ and $(Y_n)$ both start at
$x$. Recall that we write $\widehat{\eta}\equiv1-\eta$. Then by
duality, we deduce from (\ref{eqaux2}) that
%
%e60 #&#
%e61 #&#
%e62 #&#
%e63 #&#
\begin{eqnarray}\label{eqDaux}
\frac{1}{k}\Bbb E_{\mu_u}[D(x,\xi_s)]&=&-(b+c)\int\mu_u(d\eta
)\Bbb E[\widehat{\eta}(B^{X_1}_s)\eta(B^{Y_1}_s)]\nonumber\\
&&{}+b\int\mu_u(d\eta)\Bbb E[\widehat{\eta}(B^{X_1}_s)\eta
(B^{Y_1}_s)\eta(B^{Y_2}_s)]\nonumber\\
&&{}+b\int\mu_u(d\eta)\Bbb E[\eta(B^{X_1}_s)\widehat{\eta
}(B^{Y_1}_s)\widehat{\eta}(B^{Y_2}_s)]\nonumber\\
%=&-(b+c)\langle\dxi(X_1)\xi(Y_1)\rangle_{s,u}+b\langle%\dxi(e_1)
%=&-(b+c)\langle\xi(e_2)\rangle_{s,u}+(b+c)\langle%\xi(e_1)\xi(e_2)
%&+b\langle%\xi(e_1)-\xi(e_1)\xi(e_2)-\xi(e_1)\xi(e_3)+\xi(e_2)
&=&-c\int\mu_u(d\eta)\E[\eta(B^{Y_1}_s)]+c\int\mu_u(d\eta)\E[
\eta(B^{X_1}_s)
\eta(B^{Y_1}_s)]\\
&&{}+b\int\mu_u(d\eta)\E[\eta(B^{Y_1}_s)\eta(B^{Y_2}_s)]\nonumber\\
&&{}-b\int\mu_u(d\eta)\E[\eta(B^{X_1}_s)\eta(B^{Y_2}_s)]\nonumber\\
&&{}+b\int\mu_u(d\eta)\Bbb E[\eta(B^{X_1}_s)-\eta
(B^{Y_1}_s)].\nonumber
\end{eqnarray}
For clarity, let us write from now on $\mathbf P_\rho$ and $\mathbf
E_\rho$ for the probability measure and the expectation, respectively,
under which the common initial position of $(X_n)$ and $(Y_n)$ is
distributed as $\rho$. Recall that $\overline{D}(\xi)$ is the $\pi
$-expectation of $x\lmt D(x,\xi)$ defined by (\ref{eqHbar}). Write
$M_{x,y}=\inf\{t\in\Bbb R_+;B^x_t=B^y_t\}$ for the first meeting time
of the random walks $B^x$ and $B^y$, so $B^x$ and $B^y$ coincide after
$M_{x,y}$.
Then from (\ref{eqDaux}), the spatial homogeneity of the Bernoulli
product measures implies that
\begin{eqnarray*}
\frac{1}{k}\Bbb E_{\mu_u}[\overline{D}(\xi_s)]%=&-c
&=&-cu+c[u\mathbf P_\pi(M_{X_1,Y_1}\leq s)+u^2\mathbf P_\pi
(M_{X_1,Y_1}>s)]\\
&&{}+b[u\mathbf P_\pi(M_{Y_1,Y_2}\leq s)+u^2\mathbf P_{\pi
}(M_{Y_1,Y_2}>s)]\\
&&{}-b[u\mathbf P_\pi(M_{X_1,Y_2}\leq s)
+u^2\mathbf P_\pi(M_{X_1,Y_2}>s)]
\\
&=&-cu(1-u)\mathbf P_\pi(M_{X_1,Y_1}>s)
-bu(1-u)\mathbf P_\pi(M_{Y_1,Y_2}>s)\\
&&{}+bu(1-u)\mathbf P_\pi(M_{X_1,Y_2}>s).%\label{eqB2}
\end{eqnarray*}
To obtain $\B I$, we integrate both sides of the foregoing equality
with respect to $s$ over $\Bbb R_+$. This gives
%
%e64 #&#
\begin{equation}
\quad\B I(u)=ku(1-u)(-c\mathbf E_\pi[M_{X_1,Y_1}]-b\mathbf E_\pi
[M_{Y_1,Y_2}]+b\mathbf E_\pi[M_{X_1,Y_2}]).\label{eqaux8}
\end{equation}

%The formula (\ref{eqaux8}) has an interesting property to observe.
%%Informally speaking, we deduce from the Markov property and the time
%%reversibility, the expecation $\mathbf E_\pi[M_{X_1,Y_1}]$ is closely
%%related the expected length of certain path loop of length
%%$2M^+_{X_0,X_0}$. Here, $M^+_{X_0,X_0}$ is the first time $\geq1$
%that %two copies of the random walk starting at $X_0$ meet. A similar
%%interpretation applies to $\mathbf E_\pi[M_{X_1,Y_2}]$ and $\mathbf
%%E_\pi[M_{Y_1,Y_2}]$. Hence, roughly speaking, the Bernoulli transform
%%(and hence the fixation probability) is determined by a linear
%%combination of the expected lengths of some loops at $X_0$. \qed

We now turn to a simple identity between first meeting times and first
hitting times. Let $T_y=\inf\{n\geq0;X_n=y\}$, the first hitting time
of $y$ by $(X_n)$. Observe that for the random walks on (connected)
regular graphs, the invariant distribution is uniform and $\mathbf
E_x[T_y]=\mathbf E_y[T_x]$ for any $x,y$. Hence, the proof of
Proposition 14.5 in~\cite{AFIPS} implies
%
%e65 #&#
\begin{equation}
\mathbf E[M_{x,y}]=\tfrac{1}{2}\mathbf E_x[T_y], \qquad x,y\in\mathsf
V.\label{eqMT}
\end{equation}
%
%disregard of the stronger assumption therein that the underlying
%Markov chain is \textit{symmetric}.
Write
\[
f(x,y):=\mathbf E_x[T_y]=\mathbf E_y[T_x], \qquad x,y\in\mathsf V,
\]
where $\mathbf E_x=\mathbf E_{\delta_x}$.
%Using (\ref{eqMT}), we rewrite (\ref{eqaux8}) as
%[-cf(X_1,Y_1)+bf(X_1,Y_2)-bf(Y_1,Y_2)].\label{eqaux9}
%%u(1-u)(-c\Bbb E_{}[T_{e_2}]+b\Bbb E_{e_1}[T_{e_3}]-b\Bbb
%%E_{e_2}[T_{e_3}]).

%le12 #&#
\begin{lem}\label{lemct}
For any $z\in\mathsf V$,
%
%e66 #&#
%e67 #&#
%e68 #&#
\begin{eqnarray}
\mathbf E_z[f(X_0,X_1)]&=&\mathbf E_z[f(Y_1,Y_2)]=N-1,\label
{eqtime1}\\
\mathbf E_z[f(X_1,Y_1)]&=&\mathbf E_z[f(Y_0,Y_2)]=N-2,\label
{eqtime2}\\
\mathbf E_{z}[f(X_1,Y_2)]&=&\biggl(1+\frac{1}{k}\biggr)(N-1)+\frac
{1}{k}-2.\label{eqtime3}
\end{eqnarray}
\end{lem}
\begin{pf}
The proof of the equality $\mathbf E_z[f(X_0,X_1)]=N-1$ can be found in
Chapter~3 of~\cite{AFIPS} or~\cite{LRWG}. We restate its short proof
here for the convenience of readers. Let $T^+_x=\inf\{n\geq1;X_n=x\}$
denote the first return time to $x$. A~standard result of Markov chains says
$\mathbf E_x[T_x^+]=\pi(x)^{-1}=N$ for any $x$. The equalities in~(\ref{eqtime1}) now follow from the
Markov property.\vadjust{\goodbreak}

Next, we prove (\ref{eqtime2}). By (\ref{eqtime1}) and the symmetry
of $f$, we have
\begin{eqnarray*}
N-1&=&\mathbf E_z[f(X_0,X_1)]=\sum_{x\sim z}\frac{1}{k}\mathbf
E_z[T_x]\\[-3pt]
&=&\sum_{x\sim z}\sum_{y\sim z}\frac{1}{k^2}(\mathbf
E_y[T_x]+1)=\mathbf E_z[f(Y_1,X_1)]+1,
\end{eqnarray*}
so $\mathbf E_z[f(X_1,Y_1)]=N-2$. Here, our summation notation $\sum
_{x\sim z}$ means summing over indices $x$ with $z$ fixed, and the same
convention holds in the proof of (\ref{eqtime3}) and Section~\ref
{secuvmp} below.
A similar application of the Markov property to the coordinate $Y_1$ in
$\mathbf E_z[f(Y_0,Y_1)]$ gives $\mathbf E_z[f(Y_0,Y_2)]=N-2$. This
proves (\ref{eqtime2}).

Finally, we need to prove (\ref{eqtime3}). We use (\ref{eqtime1})
and (\ref{eqtime2}) to get
%
%e69 #&#
%e70 #&#
\begin{eqnarray}
\mathbf E_z[f(X_0,X_1)]&=&1+\mathbf E_z[f(X_1,Y_1)]\nonumber\\[-3pt]
&=&1+\sum_{x\sim z}\sum_{\stackrel{\scriptstyle y\sim z}{y\neq
x}}\frac{1}{k^2}\mathbf E_x[T_y]\nonumber\\[-3pt]
&=&1+\sum_{x\sim z}\mathop{\sum_{ y\sim z}}_{y\neq x}
\frac{1}{k^2}\biggl(\sum_{w\sim y}\frac{1}{k}\mathbf
E_x[T_w]+1\biggr)\label{eqauxCT1}\\[-3pt]
&=&1+\sum_{x\sim z}\mathop{\sum_{y\sim z}}_{y\neq
x}\frac{1}{k^2}+\sum_{x\sim z}\sum_{y\sim z}\sum_{w\sim y}\frac
{1}{k^3}\mathbf E_x[T_w]-\sum_{x\sim z}\sum_{w\sim x}\frac
{1}{k^3}\mathbf E_x[T_w]\nonumber\\[-3pt]%\label{eqauxCT2}\\
&=&2-\frac{1}{k}+\mathbf E_z[f(X_1,Y_2)]-\frac{1}{k}\mathbf
E_z[f(X_1,X_0)].\label{eqauxCT3}
\end{eqnarray}
Here, in (\ref{eqauxCT1}) we use the symmetry of $f$, and the last
equality follows from (\ref{eqtime1}). A rearrangement of both sides
of (\ref{eqauxCT3}) and an application of (\ref{eqtime1}) then lead
to~(\ref{eqtime3}), and the proof is complete.\vspace*{-2pt}
\end{pf}

%Similar to the proof of Lemma~\ref{lemct}, one can show that $\Bbb%E_
%transform $\B I$ can be written as

Apply Lemma~\ref{lemct} and (\ref{eqMT}) to (\ref{eqaux8}), and
we obtain the following result.\vspace*{-2pt}

%pr13 #&#
\begin{prop}\label{propBTdbu}
For any $u\in[0,1]$,
%
%e71 #&#
\begin{equation}\label{eqB3}
\B I(u)=\frac{ku(1-u)}{2}\biggl[\biggl(\frac{b}{k}-c\biggr)
(N-2)
+b\biggl(\frac{2}{k}-2\biggr)\biggr].\vspace*{-2pt}
\end{equation}
\end{prop}

Finally, since $\B I(u)$ takes the form (\ref{eqBfSpec}), we may
apply (\ref{eqBIu}) and Proposition~\ref{propBTdbu} to obtain the
explicit formula for the coefficient of $w$ in (\ref{ineqcrt}),
subject to each initial distribution $\mathbf u_n$. This proves our
assertion in Theorem~\ref{thmi1}(1).\vspace*{-2pt}

%s4.2 #&#
\subsection{\texorpdfstring{Proof of Theorem \protect\ref{thmi1}(2)}{Proof of Theorem 1(2)}}\label{subseciuT}
The proof of Theorem~\ref{thmi1}(2) follows from almost the same argument for
Theorem~\ref{thmi1}(1) except for more complicated arithmetic. For this reason,
we will only point out the main steps,\vadjust{\goodbreak} leaving the detailed arithmetic
to the interested readers. %We will only work under the special payoff
%matrix~\ref{eqPi}, leaving the general case to the reader.
In the following, we continue to use the notation for the random walks
in the proof of Theorem~\ref{thmi1}(1).

Fix $x\in\mathsf V$ and assume the chains $(X_n)$ and $(Y_n)$ both
start at $x$.
By Proposition~\ref{propiu}, we have
\begin{eqnarray*}
&&\frac{1}{k}\Bbb E_{\mu_u}[D(x,\xi_s)]\\[-2pt]
&&\qquad= \Bbb E_{\mu_u}[(b-c)\dxi_s(x)f_{11}(x,\xi_s)-c\dxi
_s(x)f_{10}(x,
\xi_s)-b\xi_s(x)f_{01}(x,\xi_s)]\\[-2pt]
&&\qquad\quad  {}-\frac{k}{k+1}\Bbb E_{\mu_u}\bigl[
\bigl((b-c)f_{11}(x,\xi_s)-cf_{10}(x,\xi_s)+bf_{01}(x,\xi_s)
\bigr)\\[-2pt]
&&\hspace*{106pt}\qquad\quad {}\times     \bigl(\dxi_s(x)f_1(x,\xi
_s)-\xi_s(x)f_0(x,\xi_s)\bigr)\bigr]\\[-2pt]
&&\qquad\quad{}-\frac{1}{k+1}\Bbb E_{\mu_u}\bigl[b\dxi_s(x)f_1^2(x,\xi_s)\\[-2pt]
&&\hspace*{92pt}{}-\xi
_s(x)f_0(x,\xi_s)\bigl(
(b-c)f_1(x,\xi_s)-cf_0(x,\xi_s)\bigr)\bigr]\\[-2pt]
&&\qquad=\int\mu_u(d\eta)\biggl(\E[ (b-c)\widehat{\eta}(B^x_s)
\eta(B^{Y_1}_s)\eta(B^{Y_2}_s)-c\widehat{\eta}(B^x_s)\eta
(B^{Y_1}_s)\widehat{\eta}(B^{Y_2}_s)\\[-2pt]
&&\hspace*{184pt}\qquad\quad{}      -b\eta(B^x_s)\widehat{\eta
}(B^{Y_1}_s)\eta(B^{Y_2}_s)]\\[-2pt]
&&\hspace*{48pt}\qquad\quad{}-\frac{k}{k+1}\Bbb E\bigl[ \bigl((b-c)\eta(B^{Y_1}_s)\eta
(B^{Y_2}_s)-c\eta(B^{Y_1}_s)\widehat{\eta}
(B^{Y_2}_s)\\[-2pt]
&&\hspace*{206pt}\qquad{}+b\widehat{\eta}(B^{Y_1}_s)\eta(B^{Y_2}_s)
\bigr)\\[-2pt]
&&\hspace*{198pt}\qquad{}\times\bigl(\eta
(B^{X_1}_s)-\eta(B^x_s)\bigr)\bigr]\\[-2pt]
&&\hspace*{48pt}\qquad\quad{}-\frac{1}{k+1}\Bbb E[ b\widehat{\eta}(B^x_s)\eta
(B^{X_1}_s)\eta(B^{Y_1}_s)\\[-2pt]
&&\hspace*{128pt}{}-(b-c)\eta(B^x_s)\widehat{\eta}
(B^{X_1}_s)\eta(B^{Y_1}_s)\\[-2pt]
&&\hspace*{216pt}{}      +c\eta(B^x_s)\widehat{\eta}(B^{X_1}_s)
\widehat{\eta}(B^{Y_1}_s)]\biggr),
\end{eqnarray*}
where the last equality follows again from duality. The last equality gives
\begin{eqnarray*}
&&\frac{1}{k}\Bbb E_{\mu_u}[D(x,\xi_s)]\\[-2pt]
&&\qquad=\int\mu_u(d\eta)\Bbb
E\biggl[ b\eta(B^{Y_1}_s)\eta(B^{Y_2}_s)+
\frac{c+b}{k+1}\eta(B^x_s)\eta(B^{Y_1}_s)\\[-2pt]
&&\hspace*{66pt}\qquad{}-c\eta(B^{Y_1}_s)-\frac{b}{k+1}\eta(B^x_s)\eta(B^{Y_2}_s)
+\frac{kc-b}{k+1}
\eta(B^{X_1}_s)\eta(B^{Y_1}_s)\\[-2pt]
&&\hspace*{66pt}\qquad{}-\frac{kb}{k+1}\eta(B^{X_1}_s)\eta(B^{Y_2}_s)+\frac{c}{k+1}\eta
(B^x_s)\eta(B^{X_1}_s)
\\[-2pt]
&&\hspace*{256pt}\qquad{}-\frac{c}{k+1}\eta(B^x_s)\biggr].
\end{eqnarray*}
Recall that $X_1\stackrel{(\rm d)}{=}Y_1$. Hence, by the definition of
$\overline{D}$ and $\mathbf P_\pi$, the foregoing implies that
\begin{eqnarray*}
\frac{1}{k}\Bbb E_{\mu_u}[\overline{D}(\xi_s)]%=&(b-c)
%[u-u(1-u)\mathbf P_\pi(M_{Y_1,Y_2}>s)]+\frac{(k+2)c+b}{k+1}
%[u-u(1-u)\mathbf P_\pi(M_{x,Y_1}>s)]\\
%&-cu-\frac{b}{k+1}[u-u(1-u)\mathbf%P_\pi(M_{x,Y_2}>s)]+
%&-
&=&-bu(1-u)\mathbf P_\pi(M_{Y_1,Y_2}>s)\\[-2pt]
&&{}-\frac{2c+b}{k+1}u(1-u)\mathbf P_\pi(M_{X_0,X_1}>s)\\[-2pt]
&&{}+\frac{b}{k+1}u(1-u)\mathbf P_\pi(M_{Y_0,Y_2}>s)\\[-2pt]
&&{}-\frac{kc-b}{k+1}u(1-u)\mathbf P_\pi(M_{X_1,Y_1}>s)\\[-2pt]
&&{}+
\frac{kb}{k+1}u(1-u)\mathbf P_\pi(M_{X_1,Y_2}>s).
\end{eqnarray*}
Again, we integrate both sides of the foregoing display with respect to
$s$ and then apply (\ref{eqMT}) and Lemma~\ref{lemct} for the
explicit form of $\B I$. The result is given by the following proposition.

%&+(b-\frac{ck}{k+1}-\frac{bk}{k+1}+\frac{b}{k+1})\mathbf%E_
%+\frac{kb}{k+1}\mathbf E_\pi[M_{X_1,Y_2}]]
%Using Lemma~\ref{lemct}, we get the following.

%pr14 #&#
\begin{prop}\label{propBTiu}
For any $u\in[0,1]$,
\[
\B I(u)
=\frac{k(k+2)u(1-u)}{2(k+1)}\biggl[\biggl(\frac{b}{(k+2)}-c
\biggr)(N-1)-\frac{
(2k+1)b-ck}{k+2}\biggr].
\]
\end{prop}

Our assertion for Theorem~\ref{thmi1}(2) now follows from an
application of Proposition~\ref{propBTiu} similar to that of
Proposition~\ref{propBTdbu} for Theorem~\ref{thmi1}(1). The proof is now complete.

%s5 #&#
\section{\texorpdfstring{Proofs of Propositions~\protect\ref{propdbu} and~\protect\ref{propiu}}
{Proofs of Propositions 2.1 and 2.2}}\label{secuvmp}
%s5.1 #&#
\subsection{\texorpdfstring{Proof of Proposition~\protect\ref{propdbu}}{Proof of Proposition 2.1}}
Suppose that $\xi\in\{1,0\}^{\mathsf V}$ is the present configuration
on the graph. Let $n_i(x)=n_i(x,\xi)$ be the number of neighboring $i$
players for an individual located at vertex $x$ for $i=1,0$. Let $w\in
[0,1]$ denote the intensity of selection. By definition, the fitness
$\rho_i(x)=\rho_i(x,\xi)$ of an $i$-player located at $x$ is given by
%
%e72 #&#
\begin{equation}\label{deffitness}
\rho_i(x)=(1-w)+w\left[
\matrix{
\Pi_{i1} & \Pi_{i0} }
\right]\left[
\matrix{
n_1(x) \vspace*{2pt}\cr
n_0(x) }
\right]=(1-w)+w\bsb\Pi_i\bsm n(x).% \mbox{ if $\xi(x)=i$}.
\end{equation}
Here, $\bsb\Pi_i$ is the payoff row of an $i$-player of the matrix
$\bsb\Pi$ and $\bsm n(x)$ is the column vector $[{n_{1}(x)\enskip
n_{0}(x)}]^\top$. Hence, there exists $w_0>0$ depending only on
$k$ and $\bsb\Pi$ such that $\rho_i>0$ for every $w\in[0,w_0]$ (see
\cite{CDPVMP}).

The game with the death--birth updating under weak selection defines a
Markov chain with transition probabilities $P^w$ taking the form (\ref
{eqcw}) and
%
%e73 #&#
%e74 #&#
\begin{eqnarray}
c^w(x,\xi)&=&r_{1-\xi(x)}(x,\xi)\geq0,\label{egcw}\\[-2pt]
r_i(x,\xi)&=&\frac{\sum_{y\sim x}\rho_i(y)\1_{\xi(y)=i}}{\sum
_{y\sim x}[\rho_1(y)\xi(y)+\rho_0(y)\dxi(y)]}.
\label{egri}
\end{eqnarray}
It has been shown in Section 1.4 of~\cite{CDPVMP} that the rates
$c^w$ define voter model perturbations satisfying (\ref{A1}) and (\ref
{A2}). Moreover, $\lambda=1$ and the functions $h_i$ in the expansion
(\ref{A2}) are given by (\ref{eqdbh}). %and
%r_i(x,\xi)=f_i(x,\xi)+wh_i(x,\xi) +w^2g_{w,i}(x,\xi)
%where
Plainly, $c^w(x,\mathbf1)\equiv r_0(x,\mathbf1)\equiv0$ and
$c^w(x,\mathbf0)\equiv r_1(x,\mathbf0)\equiv0$. Hence, (\ref{A3})
is also satisfied.

It remains to check that (\ref{A4}) is satisfied. Since (\ref{A4}) is
satisfied when $w=0$, it is enough to show that
%
%e75 #&#
\begin{equation}
P(\xi,\xi^x)>0\quad\Longleftrightarrow\quad P^w(\xi,\xi^x)>0\label{eqeq}
\end{equation}
for any $\xi\neq\mathbf1,\mathbf0$ and any $x$. However, this is
immediate from (\ref{egri}) if we notice that $\rho_i(\cdot)$ and
the constant function $1$, both regarded as measures on $\mathsf V$ in
the natural way, are equivalent.
%since $\rho_i>0$ and in the case $w=0$ all the $\rho_i$ in (
Our proof of Proposition~\ref{propdbu} is complete.
% Suppose that $\xi(x)=0$. If $P(\xi,\xi^x)=\frac{1}{N}f_1(x,\xi)>0$,
%then $x$ has a $1$-neighbor. Since $\rho_i>0$, we must have $r_1(x,
%$x$ must have a $1$-neighbor by (\ref{egri}) because $r_1(x,\xi)$ is
%the proportion of the fitness of $1$ from the $1$-neighbors of $x$.
%Hence $P(\xi,\xi^x)>0$. A similar argument applies when $\xi(x)=1$. We
%also need to claim that $P(\xi,\xi)>0$ if and only if $P^w(\xi,
%$\xi\neq\mathbf0,\mathbf1$. To see this, note we can always choose
%a vertex $x$ which is not surrounded by the strategy $1-\xi(x)$. The
%probability $r_{1-\xi(x)}(x,\xi)$ must be strictly less than $1$ by (
%positive probability the individual at the vertex $x$ is chosen to die
%and the update of $x$ gives the same strategy as before, so $P^w(\xi,

%re15 #&#
\begin{rmk}\label{rmkegfs2}
Suppose now that payoff is given by a general $2\times2$ payoff matrix
$\bsb\Pi^*=(\Pi^*_{ij})_{i,j=1,0}$ subject only to the
``equal-gains-from-switching'' condition~(\ref{eqegfs}). Let us
explain how to reduce the games with payoff matrix $\bsb\Pi^*$ to the
games with payoff matrix $\bsb\Pi^a$ under weak selection, where $\bsb
\Pi^a$ is defined by~(\ref{eqPiadj}).

In this case,
payoffs of players are as described in Remark~\ref{rmkegfs1}, and
fitness is given by
%
%e76 #&#
\begin{equation}
\rho_i^{\bsb\Pi^*}(x)=(1-w)+w\bsb\Pi^*_i\bsm n(x), \qquad x\in\mathsf
V.\label{eqrhostar}
\end{equation}
Here again, $\bsb\Pi^*_i$ is the payoff row of an $i$-player.
We put the superscript $\bsb\Pi^*$ (only in this remark) to emphasize
the dependence on the underlying payoff matrix $\bsb\Pi^*$, so, in
particular, the previously defined fitness $\rho$ in (\ref
{deffitness}) is equal to $\rho^{\bsb\Pi}$.

Suppose that the graph is $k$-regular.
The transition probabilities under the death--birth updating with payoff
matrix $\bsb\Pi^*$ are defined in the same way as before through~(\ref
{egcw}) and (\ref{egri}) with $\rho$ replaced by $\rho^{\bsb\Pi^*}$.
Note that $n_1(x)+n_0(x)\equiv k$.
Then for all small $w$
\begin{eqnarray*}
\frac{1}{1-(1-k\Pi^*_{00})w}\rho_i^{\bsb\Pi^*}(x)&=&1+\frac
{w}{1-(1-k\Pi^*_{00})w}\bsb\Pi^a_i\bsm n(x)\\
&=&1+\frac{w^a}{1-w^a}\bsb\Pi^a_i\bsm n(x)\\
&=&\frac{1}{1-w^a}\rho_i^{\bsb\Pi^a}(x)
\end{eqnarray*}
for some $w^a$. Here, $w$ and $w^a$ are defined continuously in terms
of each other by
\[
w^a=\frac{w}{1+k\Pi^*_{00}w}\quad \mbox{and} \quad w=\frac
{w^a}{1-k\Pi^*_{00}w^a},
\]
so $\lim_{w^a\to0}w=\lim_{w\to0}w^a=0$. Consequently, by (\ref
{egcw}) and (\ref{egri}), the foregoing display implies that the
death--birth updating with payoff matrix $\bsb\Pi^*$ and intensity of
selection $w$ is ``equivalent'' to the death--birth updating with payoff
matrix $\bsb\Pi^a$ and intensity of selection $w^a$, whenever $w^a$ or
$w$ is small. Here, ``equivalent'' means equality of transition probabilities.

A similar reduction applies to the imitation updating by using its
formal definition described in the next subsection, and we omit the
details.
\end{rmk}

%s5.2 #&#
\subsection{\texorpdfstring{Proof of Proposition~\protect\ref{propiu}}{Proof of Proposition 2.2}}
Under the imitation updating, the Markov chain of configurations has
transition probabilities given by
%
%e77 #&#
%e78 #&#
\begin{eqnarray}\label{eqPw}
P^w(\xi,\xi^x)&=&\frac{1}{N}\,d^w(x,\xi),
\nonumber
\\[-8pt]
\\[-8pt]
\nonumber
P^w(\xi,\xi)&=&1-\frac{1}{N}\sum_{x}\,d^w(x,\xi),
\end{eqnarray}
where
%
%e79 #&#
%e80 #&#
\begin{eqnarray}
d^w(x,\xi)&=&s_{1-\xi(x)}(x,\xi),\label{egdw}\\
s_i(x,\xi)&=&
\frac{\sum_{y\sim x} \rho_{i}(y)\1_{\xi(y)=i}}{\sum_{y\sim x}
[\rho_1(y)\xi(y)+ \rho_0(y)\dxi(y)]+\rho_{1-i}(x)}\label{egri2}
\end{eqnarray}
and the fitness $\rho_i$ are defined as before by (\ref
{deffitness}). We assume again that the intensity of selection $w$ is
small such that $\rho_i>0$.
To simplify notation, let us set the column vectors
\begin{eqnarray*}
\bsm f(x)&=&[\matrix{f_1(x) & f_0(x)}]^\top,\\
\bsm f_{i\bullet}(x)&=&[\matrix{f_{i1}(x)& f_{i0}(x)}]^\top,\\
\bsm n_{i\bullet}(x)&=&[\matrix{n_{i1}(x) & n_{i0}(x)}]^\top,
\end{eqnarray*}
where the functions $f_i$ and $f_{ij}$ are defined by (\ref{eqf}).
By (\ref{deffitness}) and (\ref{egri2}), we have
\begin{eqnarray*}
s_i(x,\xi)&=&\frac{(1-w)n_i(x)+w\bsb\Pi_i\bsm n_{i\bullet
}(x)}{(1-w)(k+1)+w\sum_{j=0}^1\bsb\Pi_j\bsm n_{j\bullet}(x)+w\bsb\Pi
_{1-i}\bsm n (x)}\\
&=&\frac{(1-w){k}/{(k+1)}f_i(x)+{k^2}/{(k+1)}w\bsb\Pi_i\bsm
f_{i\bullet}(x)}{(1-w)+w{k^2}/{(k+1)}\sum_{j=0}^1\bsb\Pi_j\bsm
f_{j\bullet}(x)+w{k}/{(k+1)}\bsb\Pi_{1-i}\bsm f(x)}\\
&=&\frac{{k}/{(k+1)}f_i(x)+w({k^2}/{(k+1)}\bsb\Pi_i\bsm
f_{i\bullet}(x)-{k}/{(k+1)}f_i(x))}{1+w(
{k^2}/{(k+1)}\sum_{j=0}^1\bsb\Pi_j\bsm f_{j\bullet}(x)+
{k}/{(k+1)}\bsb\Pi_{1-i}\bsm f(x)-1)}.
%=&\frac{\frac{k}{k+1}f_0(x)+w\theta_0(x)}{1+w\phi_0(x)}.
\end{eqnarray*}
Note that the functions $f_i$ and $f_{ij}$ are uniformly bounded. Apply
Taylor's expansion in $w$ at $0$ to the right-hand side of the
foregoing display. We deduce from~(\ref{egdw}) that the transition
probabilities (\ref{eqPw}) takes the form (\ref{eqcw}) with
$\lambda=\frac{k}{k+1}$ and the rates~$c^w$ satisfying (\ref{A1})
and (\ref{A2}) for some small $w_0$. Moreover, the functions $h_i$ are
given by
\begin{eqnarray*}
h_i&=&(k\bsb\Pi_i\bsm f_{i\bullet}-f_i)-f_i\Biggl(\frac
{k^2}{k+1}\sum_{j=0}^1\bsb\Pi_j\bsm f_{j\bullet}+\frac{k}{k+1}\bsb
\Pi_{1-i}\bsm f-1\Biggr)\\
%=&k\bs\Pi_i\bs f_{i\bullet}-f_i-\frac{k^2}{k+1}\sum_{j=0}^1f_i\bs\Pi_j
&=&k\bsb\Pi_i\bsm f_{i\bullet}-\frac{k^2}{k+1}f_i \Biggl(\sum
_{j=0}^1\bsb\Pi_j\bsm f_{j\bullet}\Biggr)-\frac{k}{k+1}f_i\bsb\Pi
_{1-i}\bsm f.
\end{eqnarray*}
By the definition of $\bsb\Pi$ in (\ref{eqPi}), we get (\ref{eqiuh}).
The verifications of (\ref{A3}) and (\ref{A4}) follow from similar
arguments for those of (\ref{A3}) and (\ref{A4}) under the
death--birth updating, respectively. This completes the proof of
Proposition~\ref{propiu}.

\section*{Acknowledgments} I am grateful to my advisor, Professor Ed
Perkins, for inspiring discussions and to Christoph Hauert and the
anonymous referees for valuable comments.

% imsref loaded by akundreckaite, 2012-06-05 12:35:42
%

%suskaldyti doi

\printaddresses

\end{document}